\documentclass[reqno, 12pt]{amsart}

\usepackage{hyperref}

\usepackage{enumerate}

\usepackage{amssymb, amsmath}
\usepackage{epsfig}
\newcommand\bx{{\mathbf x}}
\newcommand\by{{\mathbf y}}
\newcommand\bz{{\mathbf z}}
\newcommand\bq{{\mathbf q}}
\newcommand\bp{{\mathbf p}}
\newcommand\be{{\mathbf e}}

\newcommand\bk{{\mathbf k}}

\newcommand\br{{\mathbf r}}
\newcommand\bX{{\mathbf X}}
\newcommand\bK{{\mathbf K}}

\newcommand\C{{\mathbb C}}

\newcommand\R{{\mathbb R}}
\newcommand\T{{\mathbb T}}
\newcommand\Z{{\mathbb Z}}
\newcommand\ve{\varepsilon}

\newcommand\bpsi{\boldsymbol{\psi}}

\newcommand\tom{\tilde{\omega}}

\newcommand\mcE{\mathcal E}
\newcommand\mcW{\mathcal W}
\newcommand\mcY{\mathcal Y}

\newtheorem{assumption}{Assumption}
\newtheorem{theorem}[assumption]{Theorem}
\newtheorem{remark}[assumption]{Remark}
\newtheorem{proposition}[assumption]{Proposition}
\newtheorem{lemma}[assumption]{Lemma}

\begin{document}
\title[Energy Transport]{
Energy Transport in Stochastically Perturbed Lattice Dynamics}

\author{Giada Basile}\email{basile@wias-berlin.de}
\address{WIAS, Mohrenstr. 39, 10117 Berlin, Germany}
\author{Stefano Olla}
\email{olla@ceremade.dauphine.fr}
\address{Ceremade, UMR-CNRS 7534, Universit\'e de Paris Dauphine,
75775 Paris Cedex 16, France.
}
\author{Herbert Spohn}
 \email{spohn@ma.tum.de}
\address{Zentrum Mathematik, TU M\"{u}nchen, D-85747 Garching, Germany}

\date{\today}

\begin{abstract}
We consider lattice dynamics with a small stochastic perturbation of
order $\varepsilon$ and prove that for a space-time scale of order
$\varepsilon^{-1}$ the local spectral density (Wigner function)
evolves according to a linear 
transport equation describing inelastic collisions. For an energy
and momentum conserving chain the transport equation predicts a slow
decay, as $1/\sqrt t$, for the energy current correlation in
equilibrium. This is in agreement with previous studies using a
different method.
\end{abstract}

\thanks{S.O. research was supported by 
French ANR LHMSHE, n.BLAN07-2184264.}
\maketitle

\section{Introduction}
\label{sec:introduction}

Over the recent years there have been strong efforts to understand
energy transport in anharmonic chains. A prototypical example is the
Fermi-Pasta-Ulam chain with pinning potential and a quartic
nonlinearity. Its hamiltonian reads
\begin{equation}\label{1}
H({q},{p})=
\sum^N_{j=-N}\{\tfrac{1}{2m}p^2_j+\tfrac{1}{2}\omega^2_0 q^2_j\}+
\sum^{N-1}_{j=-N}\{\tfrac{1}{2}\alpha_1(p_{j+1}-q_j)^2+\beta(q_{j+1}-q_j)^4\}\,.
\end{equation}
Here $q_j$ is the displacement of the $j$-th particle from its
equilibrium position $j$, $p_j$ is the canonically conjugate
momentum, $m$ is the mass of the particle. 
The potential has two parts. The term $\omega^2_0 q^2_j$ confines the
$j$-particle and is therefore referred to as pinning potential. The
second part depends only on displacement differences. It consists of a
harmonic nearest neighbor term of strength $\alpha_1>0$, and the
anharmonicity $\beta( q_{j+1}-q_j)^4, \beta\ge 0$.
 Energy
transport can be studied by coupling the end particles, $q_{-N}$ and
$q_N$, to thermal reservoirs at different temperatures. An
alternative method would be to take first $N\to\infty$ in thermal
equilibrium and to monitor the spreading of the energy when
initially deposited close to the origin. If $\beta=0$, the energy
spreading is ballistic and the steady state energy flux, $j_N$, is
independent of $N$. For $\beta> 0$ and for non-zero pinning one
finds, mostly based on numerical simulations, that the current
$j_N=\mathcal{O}(1/N)$. This is called regular transport.
Generically the spreading of energy is then diffusive. On the other
hand, for $\omega_0=0$, the energy transport is anomalous,
$j_N=N^{-\alpha}$ with $\alpha$ in the range from $1/3$ to $2/5$
(see \cite{sll} for a review on the subject).

The mathematical analysis of the energy transport in the FPU-chain
is a difficult task, in particular since non-zero temperature is
required. In \cite{bo} one of the author's proposed to replace the
nonlinearity in (\ref{1}) by a stochastic exchange of
momentum between neighboring sites such that the local
 energy is conserved, with the extension to local momentum conservation
 being worked out in \cite{bborev, bbo2}.
We refer to section \ref{sec:mod1d} below for a precise 
definition. The goal of our paper is to understand in this model the
mechanism of energy transport for small noise strength, in
particular to determine conditions for regular, resp. anomalous,
transport.

Energy transport refers to large space-time scales. In case
$\beta=0$ the equations of motion for (\ref{1}) reduce to a discrete
wave equation, and semiclassical analysis provides a convenient tool. 
For the harmonic chain at non-zero temperature such a program was carried
out by Dobrushin \textit{et al.} in a series of papers \cite{dobru,
  dobru1, dobru2}. The local spectral density $W$, in physics
parlance the Wigner function, is governed by the linear transport
equation
\begin{equation}\label{2}
\frac{\partial}{\partial t} W(r,k,t)
+\frac{1}{2\pi}\omega'(k)\frac{\partial}{\partial r}W(r,k,t)=0\,.
\end{equation}
$W$ depends on position $r\in \mathbb{R}$ along the chain, the wave
number $k\in [-\frac{1}{2},\frac{1}{2}]$, and time $t$. In our model
the dispersion relation $\omega$ is computed from
\begin{equation}\label{3}
\omega(k)^2=\omega^2_0+\alpha_1\big(1-\cos(2\pi k)\big)\,.
\end{equation}
The local energy density at time $t$ is defined through
\begin{equation}\label{4}
e(r,t)=\int^{1/2}_{-1/2}  W(r,k,t) dk
\end{equation}
in our units. Thereby (\ref{2}) provides a quantitative description for
the energy transport in the harmonic chain.

In \cite{ds} the analysis of Dobrushin et al. is extended to a wider
class of harmonic lattices, including higher dimensions. Mielke
\cite{Mi} covers at great generality the case of initial conditions of
finite energy, see also \cite{LST}.

The starting point of our investigation is the observation that the
local spectral density should also be the appropriate quantity when
stochastically perturbing the lattice dynamics, provided the strength
of the stochastic term is appropriately adjusted to the semiclassical
limit. This program indeed works out and, as 
 our main result, we will prove that the stochastic exchange gives
rise to a linear collision term in (\ref{2}) and the transport
equation is modified to 
\begin{equation}
  \label{2C}
  \frac{\partial}{\partial t} W(r,k,t)
  +\frac{1}{2\pi}\omega'(k)\frac{\partial}{\partial r}W(r,k,t)=
  \int R(k,k') \left(W(r,k',t) - W(r,k,t)\right) \; dk \,. 
\end{equation}
with an explicit transition kernel $R(k,k')$.

 We will give in detail the proofs in the one-dimensional case
with a stochastic dynamics that conserves total momentum and total
energy (as in \cite{bborev}). These proofs can be straightfowardly
extended to lattices of dimension $d\ge 2$, and to stochastic
perturbations that conserves only energy (as in \cite{bo}).  

The limit transport equation is easily analysed and is a convenient
tool to study energy transport. 
For the one dimensional model with conservation of energy and
momentum, one finds that for $\omega_0\neq 0$ the energy transport is
diffusive while for no pinning, $\omega_0=0$, the spreading is
superdiffusive, in accord with the results in \cite{bborev, bbo2}.
The superdiffusion process is governed by a Levy process of index
$\frac{3}{2}$ (see \cite{JO}).
If  $\omega_0  > 0$ one finds regular energy transport. the same holds for dimension $d\ge 3$ and for noise conserving only energy, 
in accordance with a direct approach using the Green-Kubo formula
  \cite{bborev, bbo2, bo}. 


Our paper is organized as follows. In sections \ref{sec:mod1d} to 
\ref{proof} we derive equation
(\ref{2}) including the collision term. In section
\ref{sec:extens-d-dimens}  we explain how our proof
generalizes to higher dimensional lattices. 
  We then discuss the energy
transport in  detail (section \ref{sec:homog-case:-corr}).

\section{The model (One-dimensional case)}
\label{sec:mod1d}

To develop the necessary techniques, we consider first the case of a
one-dimensi\-onal chain, as in (\ref{1}) with $\beta = 0$. 
In view of higher dimension it is convenient to label the particles by
$y\in\Z$. The chain is infinite and we allow for a general harmonic
coupling.   
The phase space is $(\R\times\R)^\Z$ and a configuration at time $t$ is
denoted by $\{q_y(t),p_y(t)\}_{y\in\Z}$.
The Hamiltonian of the system is given by
\begin{equation}
\label{def:Hamilt}
H(p,q) = \frac 12 \sum_{y\in\Z} {p_y^2} +
\frac 12 \sum_{y,y'\in\Z} \alpha(y - {y'})
q_y q_{y'},
\end{equation}
where we use units such that $m=1$.
We denote with $\hat{v}(k)$, $k\in\T = [0,1]$, the Fourier transform
of a function $v$ on $\Z$,
\begin{equation}
\hat{v}(k)=\sum_{z\in\Z}e^{-2\pi i k z} v(z),
\end{equation}
and with $\tilde f(z)$, $z\in\Z$, the inverse Fourier transform
of a function $f$ on $\T$,
\begin{equation}
\tilde{f}(z)=\int_{\T}dk \; e^{2\pi ik z} f(k).
\end{equation}
The function $\omega(k)=\sqrt{\hat{\alpha}(k)}$ is called
\emph{dispersion relation}.


We assume 
$\alpha(\cdot)$ to satisfy the following properties:
\begin{assumption}\label{alpha_assump}\hspace{4cm}
\begin{itemize}
\item\emph{(a1)}
$\alpha(y)\neq 0$ for some $y\neq 0$.
\item\emph{(a2)}
$\alpha(y)=\alpha(-y)$ for all $y\in\Z$.
\item \emph{(a3)}
There are constants $C_1,C_2>0$ such that for all $y$
$$
|\alpha(y)|\leq C_1e^{-C_2|y|}.
$$
\item \emph{(a4)}
We require either
\begin{itemize}
\item \emph{(pinning)}: $\hat{\alpha}(k)>0$ for all $k\in\T$ 
\item \emph{(no pinning)}: $\hat{\alpha}(k)>0$ $\forall k\neq 0$,
  $\hat{\alpha}(0)=0$,  $\hat{\alpha}''(0) > 0$.
\end{itemize}
\end{itemize}
\end{assumption}
%
Assumptions (a2), (a3) ensure that $\hat{\alpha}$ is a real analytic
function on $\T$. Note that the potential depends only on the
displacements differences iff $\omega(0)=0$, hence the two cases in
(a4). In the pinned case $\omega$ is
strictly positive and thus real analytic on $\T$. In the unpinned case,
the condition (a4) says that $\omega(k) = c |k|$ with $c>0$ for
small $k$, to say $\omega$ is a \emph{regular acoustic} dispersion
relation.


We consider the Hamiltonian dynamics weakly perturbed by a stochastic
noise acting only on momenta and locally preserving  momentum and
kinetic energy.
The generator of the dynamics is
\begin{equation}
L=A+\ve\gamma S
\end{equation}
with $\ve > 0$, where $A$ is the usual Hamiltonian vector field
\begin{equation}
\label{def:Agen}
\begin{split}
A = \sum_{y\in\Z} p_{y} \partial_{q_{y}} - \sum_{y,y'\in\Z} \alpha(y
- y') q_{y'}  \partial_{p_{y}},
\end{split}
\end{equation}
while $S$ is the generator of the stochastic perturbation. The
operator $S$ acts only on the momenta $\{p_y\}$ and generates
a diffusion on the surface of constant kinetic energy and constant
momentum. $S$ is defined as
\begin{equation}\label{defq:Sgen1}
S = \frac 16 \sum_{z\in\Z}(Y_z)^2,
\end{equation}
where
\begin{equation*}
\label{eq:005}
Y_z=(p_z-p_{z+1})\partial_{p_{z-1}}+(p_{z+1}-p_{z-1})\partial_{p_z} +
(p_{z-1}-p_z)\partial_{p_{z+1}}
\end{equation*}
which is a vector field tangent to the surface of constant kinetic
energy and of constant momentum for three neighbouring particles. As
a consequence energy and momentum are locally conserved which, of course, implies also the conservation of total momentum and total
energy of the system,
\begin{equation*}
S \; \sum_{y\in\Z} p_y = 0\ ,\hspace{0.4cm} S H = 0 .
\end{equation*}


The evolution of $\{p(t),q(t)\}$ is given by the following
stochastic differential equations
\begin{equation}
\label{eq:sde1}
\begin{split}
dq_y= & \;p_y \;dt ,\\
dp_y = &
-(\alpha * q)_{y}\; dt + \frac{\ve\gamma} 6
\Delta(4p_y + p_{y-1} + p_{y+1}) dt \\ &
+ \sqrt{\frac{\ve \gamma}3}
\sum_{k=-1,0, 1} \left(Y_{y+k} p_y \right) dw_{y+k}(t).
\end{split}
\end{equation}
Here $\{w_{y}(t)\}_{y\in\Z}$ are independent standard Wiener
processes and $\Delta$ is the discrete laplacian on $\Z$,
\begin{equation*}
\Delta f(z) = f(z+1) + f(z-1) -2f(z) .
\end{equation*}

To study the local spectral density it is convenient to introduce the
complex valued field $\psi:\Z\to\C$ defined as 
\begin{equation}
\label{def:psi}
\begin{split}
\psi(y,t) = \frac{1}{\sqrt{2}}\big((\tilde{\omega}\ast q)_y(t)
+ip_y(t)\big).
\end{split}
\end{equation}
Observe that $|\psi(y)|^2 =
\frac{1}{2}p^2_y+\frac{1}{2}\sum_{y'\in\Z}\alpha(y-y')q_y q_{y'}=
e_y$ is the energy of particle $y$ and conservation of total energy
is equivalent to the conservation of the $\ell_2$-norm. 
For every $t\geq 0$ the evolution of $\psi$ is given by the
stochastic differential equations,
\begin{equation}\label{evol}
\begin{split}
d\psi(y,t)= & -i(\tom\ast \psi)(y,t) dt +\frac{1}{2}\ve\gamma\beta\ast(\psi -
\psi^*)(y,t) dt\\
& + \sqrt{\frac{\ve\gamma}3} \sum_{k=-1,0, 1} \big(Y_{y+k}
\tfrac{1}{2}(\psi-\psi^*)(y,t)\big) dw_{y+k}(t) ,
\end{split}\end{equation}
where $\beta$ is  defined through
\begin{equation}
\label{beta}
(\beta * f)(z)=
\frac 16 \Delta(4f(z) + f(z-1) + f(z+1)) .
\end{equation}


\section{Wigner distribution and the Boltzmann
Phonon  Equation}
\label{wignerd}

Given a complex valued function $J$ on $\R\times\T$, we define on
$\R\times\Z$
\begin{equation}
\label{eq:5} \tilde J(x, z) = \int_{\T}dk\; e^{2\pi i k z} J(x,k)
.
\end{equation}
We also define on $\R\times\T$
\begin{equation}
\label{eq:ftc} \widehat J(p, k) = \int_{\R}dx\; e^{-2\pi i p x} J(x,k)
  .
\end{equation}
We choose a class of test-functions $J$ on $\R\times\T$ such that
$J(\cdot,k)\in \mathcal{S}(\R,\C)$ for any $k\in\T$.

Let us fix $\ve>0$. We denote by $\big<\cdot\big>_\ve$
the expectation value with respect
to a family of probability measures on
phase space which satisfies the
following properties:
\begin{enumerate}[(b1)]
\item $\big<\psi(y)\big>_\ve=0,\hspace{0.4cm}\forall y\in\Z$;
\item $\big<\psi(y')\psi(y)\big>_\ve=0,\hspace{0.4cm}\forall y, y'\in\Z$;
\item $\sup_{\ve>0}\;\ve\;\big<\|\psi\|^2_{\ell_2}\big>_\ve\leq K$
 for some $K>0$.
\end{enumerate}

Observe that, since $\big< \|\psi\|^2\big>_\ve=\big< H \big>_\ve$ is
the expectation value of the energy, we are considering states
with an energy of order $\ve^{-1}$. 

For every $\psi$ we define the associated Wigner function in the
standard way, see \cite{bpr,LP,LS,RPK}, and integrate it against the
test function $J$. The bilinear expression in $\psi$ is averaged over
$<\cdot>_\ve$. For simplicity we call the averaged Wigner function
simply Wigner distribution and denote it by $W^\ve$ with $\ve$ the
small semiclassical parameter. Thereby we arrive at the following definition
\begin{equation}\label{eq:2}
\begin{split}
\big<J,W^\ve \big>= &\frac\ve 2 \sum_{y, y'\in\Z}
\big<\psi(y')^*\psi(y)\big>_\ve \int_{\T}\; dk e^{2\pi i
k(y'- y)} J (\ve(y'+y)/2,k)^*\\
= &\frac \ve 2\sum_{y, y'\in\Z} \big<\psi(y')^*\psi(y)\big>_\ve
\tilde J (\ve(y'+ y)/2,y- y')^*,
\end{split}
\end{equation}
where $J\in \mathcal S (\R\times\T)$. By condition (b3) this
distribution is well defined, see proposition \ref{prop:bound} in
the appendix.



\begin{proposition}\label{prop:bound2}
Under the assumption
(b3) for every test function $J\in \mathcal{S}(\R\times\T)$, there
exist constants $K_1, K_2$ such that
\begin{equation}
\begin{split}
\big|\big<J,W^\ve \big>\big|\leq K_1 \int_{\R}dp \;
\sup_{k \in \T}|\hat{J}(p, k)|<\infty,\\
\big|\big<J,W^\ve \big>\big|\leq K_2 \sum_{z\in\Z} \; \sup_{x \in
\Z}|\tilde{J}(x, z)|<\infty
\end{split}
\end{equation}
for every $\ve>0$.
\end{proposition}

\begin{remark}\label{rem}
Notice that $W^\ve$ is well defined on a wider class of test functions
than $\mathcal{S}(\R^d\times\T^d)$. In particular we can take $J(x,
k) = J(k)$, a bounded real valued function on $\T$, and by
(\ref{eq:2}) we have
\begin{equation*}
\big<J,W^\ve \big> = \frac \ve 2 \int_{\T} dk \;
\big<|\hat\psi(k)|^2\big>_\ve J(k) ,
\end{equation*}
while choosing $J(x, k) = J(x)$, a bounded real valued function on
$\R$, we have
\begin{equation*}
\big<J,W^\ve \big> = \frac\ve 2 \sum_{y \in \Z}\;
\big< e_y\big>_\ve J(\ve y) .
\end{equation*}
\end{remark}

Let us start our dynamics with an initial measure satisfying
conditions
(b1), (b2), (b3).
We want to study the evolution of the Wigner
distribution $W^\ve$ on the time scale $\ve^{-1}t$, i.e. we define for
$J \in \mathcal{S}(\R\times\T)$,
\begin{equation}
\begin{split}
&\big<J,W^\ve(t) \big>\\
&\hspace{8pt}=\frac \ve 2 \sum_{y, y'\in \Z}\big<\psi(y', t/\ve)^*\psi(y,
t/\ve)\big>_\ve\int_{\T} dk \;e^{2\pi i k (y'-y)} J(\ve(y'+y)/2,k)^*.
\end{split}
\end{equation}
Observe that since the dynamics preserves the total energy, the
condition $\ve \left<\|\psi \|\right>$ $ \le K$ holds at any time and, by proposition \ref{prop:bound}, the Wigner distribution is
well defined at any time.

According to remark \ref{rem}, if we choose test functions depending
only on $k$, then we obtain the distribution of energy in $k$-space.
It turns out that in the limit as $\epsilon \to 0$, this
distribution converges to the solution of the homogeneous Boltzmann
equation, and this will be our first result. We define the distribution
$\mathcal{E}^\varepsilon (t)$ on $\T$ by
\begin{equation}
 \label{eq:k-en}
 \big< J, \mcE^\ve(t) \big> = \big<J,W^\ve(t) \big> = \frac \ve 2
 \int_{\T} dk \; \big<|\hat\psi(k, t/\ve)|^2\big>_\ve J(k)
\end{equation}
for any bounded real valued function $J$, and introduce the collision
operator $C$ acting on $\mathcal{S}(\R\times\T)$,
\begin{equation}\label{C1}
C J(x,k) = \int_{\T} dk' R(k,k') \big(J(x,k')- J(x,k) \big) .
\end{equation}
with
\begin{equation}\label{kernel1}
\begin{split}
R(k,k') = \frac{4}{3}\big(2\sin^2(2\pi k)\sin^2(\pi k')+ &
2\sin^2(2\pi k')\sin^2(\pi k)\\
&-\sin^2(2\pi k)\sin^2(2\pi k')\big).
\end{split}
\end{equation}

\begin{theorem}
 \label{hom1-th}
 Assume that the initial measure satisfies conditions
 (b1),(b2), (b3), and furthermore $ \mcE^\ve(0)$ converges to a
 positive measure $\mcE_0(dk)$ on $\T$. Then  $\mcE^\ve(t)$ converges
 to $\mcE(t, dk)$, the solution of
 \begin{equation}
   \label{eq:bpeh1}
   \partial_t \big< J, \mcE(t) \big> = \gamma\big< CJ, \mcE(t) \big>
 \end{equation}
for every bounded function $J:\T\to\R$. 

\end{theorem}

In order to prove the next theorem, the full inhomogeneous equation,
we need an additional condition on the initial distribution in the unpinned
case ($\hat\alpha(0)=0$):
\begin{enumerate}[(b4)]
\item
In the case of no pinning we require
$$
\lim_{R\to 0}\,\overline{\lim_{\ve\to 0}} \; \frac\ve 2
\int_{|k|<R}dk \; \big<|\hat{\psi}(k)|^2\big>_\ve
= 0.
$$
\end{enumerate}
Condition (b4) ensures that there is no initial concentration of
energy at wave number $k=0$. This condition can be omitted for
 a dispersion relation $\omega$ which is analytic on $\T$
(as in the pinned case).

\begin{theorem}\label{th-inhbe}
Let \emph{ Assumptions (b1-b4)} hold and assume that $W^\ve(0)$
converges to a positive measure $\mu_0(dx,dk)$. Then, for all $t\in
[0,T]$, $W^\ve(t)$ converges to a positive measure $\mu(t) =\mu(t,dx,dk)$,
which is the unique solution of the Boltzmann equation
\begin{equation}
\label{be}
\begin{split}
\partial_t \big< J, \mu(t) \big>  =
\frac{1}{2\pi}\big<
\omega'(k) \partial_{x} J, \mu(t) \big>
+ \gamma \big< C J, \mu(t) \big>
\end{split}
\end{equation}
with initial condition
$\mu(0, dx, dk) = \mu_0(dx,dk)$.
\end{theorem}
\noindent
In (\ref{be}) $\big< J, \mu(t) \big>$ denotes the linear functional
$\int_{\R\times\T} J(x,k)^*\mu(t,dx,dk)$.

Observe that the kernel $R$ of (\ref{kernel1}) is non-negative,
symmetric, and is equal to zero only if $k=0$ or $k'=0$. Moreover,
it is easy to see that
\begin{equation}\label{beta1}
\int_\T dk'\; R(k,k')=-\hat{\beta}(k) ,
\end{equation}
where $\hat{\beta}(k)$ is the Fourier transform of the function
$\beta$ defined in (\ref{beta}). Thus the  Boltzmann equation
(\ref{be}) can be interpreted as the forward equation of a Markov
process $(X(t),K(t))$ on $\R\times\T$ for the dynamics of a
particle, which in the context of lattice dynamics is called
\emph{phonon}. The phonon with momentum $k$ travels with velocity
$\omega'(k)/2 \pi$ and suffers random collisions. More precisely $K(t)$ is
an autonomous reversible jump Markov process with jump rate $R$,
while the position $X(t)$ is determined through
\begin{equation*}
X(t) = X(0) + \int_0^t \frac 1{2\pi}\omega'(K(s)) \; ds.
\end{equation*}


\section{Proof of theorems \ref{hom1-th} and \ref{th-inhbe}.}
\label{proof}

\subsection{Relative Compactness of the Wigner distribution}
Existence of the limit of the Wigner distributions $W^\ve(t)$ will be
established as in \cite{bpr}. The limit
distribution $W(t)$ is non-negative, as it is proved in
\cite{LP}, \cite{LS}, i.e.
for every $J\in\mathcal{S}(\R\times \T,\C)$ one has
\begin{equation}\label{non-negW}
\big<|J|^2,W^\ve\big>= \frac\ve 2 \big<\int_{\T}dk\; \Big|\sum_{y \in \Z}
J(\ve y, k) e^{-2\pi ik y}\psi(y)\Big|^2 \big>_\ve + \mathcal{O}(\ve) .
\end{equation}

Let us introduce the space $\mathcal{A}$ of functions $J$ on
$\R^d\times\T^d$ such that
\begin{equation}
\|J\|_\mathcal{A}=\sum_{z \in \Z}\sup_{x\in\R}|\tilde{J}(x,z)| <\infty ,
\end{equation}
where $\tilde{J}$ is defined in (\ref{eq:5}). The following lemma
shows that if the distributions $W^\ve$ are uniformly bounded in
$\mathcal{A}'$, the dual space to $\mathcal{A}$, then at every time
$t$ one can choose a sequence $\ve_j\to 0$ such that $W^{\ve_j}$
converge in the *-weak topology in $\mathcal{A}'$ to a limit
distribution $W(t)$.

\begin{lemma}\label{lemma:conv_W}
There exists a constant $C>0$ independent of $t$ such that $\forall
\ve>0$
$$
\|W^\ve(t)\|_{\mathcal{A}'}\leq C.
$$
\end{lemma}
\emph{Proof.}
For every $J\in\mathcal{A}$
\begin{equation*}
\big<J,W^\ve(t)\big> = \frac \ve 2 \sum_{y, y'\in \Z}
\big<\psi(y')\psi(y) \big>_\ve \tilde{J}(\ve (y+y')/2,y-y')^* .
\end{equation*}
Then, using Schwarz inequality and Assumption (b3),
\begin{equation*}
\left| \big<J,W^\ve(t)\big>_\ve\right|\leq
\frac\ve 2\sum_{y \in \Z}
\big<|\psi(y)|^2\big>_\ve\|J\|_\mathcal{A}\leq K \|J\|_\mathcal{A} .
\end{equation*}
\qed

\medskip


\subsection{Proof of theorem \ref{hom1-th}.}

We consider a class of test functions $J$ depending only on
$k\in\T$. In particular  we  choose $J$ real valued and bounded.
Recall the definition
\begin{equation*}\begin{split}
\big<J,\mcE^\ve \big>
=\frac \ve 2 \int_{\T} dk \;
\big<|\hat\psi(k)|^2\big>_\ve J(k) ,
\end{split}\end{equation*}
which is well defined since
$|\big<J,\mcE^\ve \big>|\leq \frac 1 2 K\sup_{k\in\T}|J(k)|$.

The evolution of the distribution $\mcE^\ve(t)$ is determined by Ito's formula, namely
\begin{equation*}\begin{split}
\partial_t \big<J,\mcE^\ve(t)\big>= \frac \ve 2 \int_{\T} dk \;
\ve^{-1}\big<L|\hat\psi(k,t/\ve)|^2\big>_\ve J(k)\ ,
\end{split}\end{equation*}
where
$$
L|\hat\psi(k)|^2=A|\hat\psi(k)|^2+\ve\gamma S|\hat\psi(k)|^2
$$
and $A$, $S$ are respectively defined in (\ref{def:Agen}), (\ref{defq:Sgen1}).
We have
\begin{equation*}
A|\hat\psi(k)|^2=[A\hat{\psi}(k)^*]\hat\psi(k)+\hat\psi(k)^*[A\hat\psi(k)],
\end{equation*}
where by direct computation
\begin{equation}
A\hat\psi(k)=\sum_{y\in\Z}e^{-2\pi i ky}A\psi(y)=-i\sum_{y\in\Z}e^{-2\pi i
 ky}(\tom * \psi)(y)
=-i\omega(k)\hat\psi(k)
\end{equation}
and thus
$A|\hat\psi(k)|^2=0$.
For the stochastic part, since $S$ is a second order operator, we have
$$
S|\hat\psi(k)|^2=[S\hat{\psi}(k)^*]\hat\psi(k)+\hat\psi(k)^*[S\hat\psi(k)]
+\frac{1}{3}\sum_{z\in\Z}[Y_z\hat\psi(k)^*][Y_z\hat\psi(k)] ,
$$
where by direct computation
\begin{equation}\begin{split}
S\hat\psi(k)=&\sum_{y\in\Z}e^{-2\pi i
 ky}S\psi(y)=\frac{1}{2}\sum_{y \in \Z}e^{-2\pi i ky}\beta
 *(\psi-\psi^*)(y)\\
=&\frac 1 2\hat{\beta}(k)(\hat\psi(k)-\hat\psi(-k)^*)
\end{split}\end{equation}
with $\beta$ defined in (\ref{beta}).
Thus
\begin{equation*}\begin{split}
&(S\hat{\psi}(k)^*)\hat\psi(k)+\hat\psi(k)^*(S\hat\psi(k))\\
&\hspace{10pt}=\hat\beta
|\hat\psi(k)|^2 -\frac 1 2 \hat\beta(\hat\psi(k)\hat\psi(-k) + \hat\psi(k)^*\hat\psi(-k)^*),
\end{split}\end{equation*}
where
\begin{equation}\label{hbeta}
\hat\beta(k)=-\frac{4}{3}\sin^2(\pi k)(1+2\cos^2(\pi k)).
\end{equation}

Finally we have to compute
$\frac{1}{3}\sum_{z\in\Z}[Y_z\hat\psi(k)^*][Y_z\hat\psi(k)]$. It
holds
\begin{equation}\label{syy0}\begin{split}
\sum_{z\in\Z}[Y_z\hat\psi(k)^*][Y_z\hat\psi(k)]=\sum_{y,y'\in\Z}
e^{2\pi i k (y'- y)}\sum_{z\in\Z} [Y_z\psi(y')^*][Y_z\psi(y)],
\end{split}\end{equation}
where
\begin{equation*}
\begin{split}
&\sum_{z\in\Z} [Y_z\psi(y')^*][Y_z\psi(y)]\\
&\hspace{10pt}=  [Y_{y+1}\psi(y')^*][Y_{y+1}\psi(y)]+
[Y_y\psi(y')^*][Y_y\psi(y)]
+[Y_{y-1}\psi(y')^*][Y_{y-1}\psi(y)].
\end{split}\end{equation*}
This expression is explicitly computed in the appendix, see eq.
(\ref{syy2}). By inserting it in (\ref{syy0}) we get
\begin{equation*}\begin{split}
&\sum_{y,y'\in\Z}
e^{2\pi i k (y'- y)}\sum_{z\in\Z} [Y_z\psi(y')^*][Y_z\psi(y)]\\&
\hspace{10pt}=
\cos(4\pi k) \sum_{y\in\Z} (2p_y p_{y+1} -
p_{y}p_{y+2} - p_y^2)\\
&\hspace{24pt}+ \cos(2 \pi k) \sum_{y\in\Z}  (2p_{y}p_{y+2} - 2p_y^2)
+\sum_{y\in\Z} (-2p_y p_{y+1} - p_{y}p_{y+2} + 3p_y^2) ,
\end{split}\end{equation*}
which is equal to
\begin{equation*}\begin{split}
\int_\T d\xi\; |\hat{p}(\xi)|^2 & \big(
\cos(4\pi k)[2\cos(2\pi\xi)-\cos(4\pi\xi)-1]\\
&+2\cos(2\pi k)[\cos(4\pi \xi)-1]
+[3-2\cos(2\pi\xi)-\cos(4\pi\xi)]
\big).
\end{split}\end{equation*}
Finally, after some trigonometric identities and using the relation
$$
|\hat p(k)|^2=\frac 1 2 \big(|\hat\psi(k)|^2+ |\hat \psi(-k)|^2 -
\hat\psi(k) \hat\psi(-k)-\hat\psi(k)^*\hat\psi(-k)^*\big),
$$
we get
\begin{equation*}\begin{split}
\frac{1}{3}\sum_{z\in\Z}[Y_z\hat\psi(k)^*][Y_z\hat\psi(k)]\hspace{6cm}\\
=\int_\T d\xi\; R(k,\xi) \big(|\hat{\psi}(\xi)|^2-\frac 12
[\hat{\psi}(\xi)\hat{\psi}(-\xi)+\hat{\psi}(\xi)^*\hat{\psi}(-\xi)^*]\big),
\end{split}\end{equation*}
where $R(k,\xi)$ is given by (\ref{kernel1}).

Since
$\int_\T d\xi\; R(k,\xi)=-\hat{\beta}(k)$,
we can write
\begin{equation*}
S|\hat{\psi}(k)|^2= C|\hat{\psi}(k)|^2-\frac 1 2
C(\hat{\psi}(k)\hat{\psi}(-k)+\hat{\psi}(k)^*\hat{\psi}(-k)^*),
\end{equation*}
where $C$ is the operator defined in (\ref{C1}), i.e.
$$
Cf(k)=\int_\T d\xi\; R(k,\xi)\big(f(\xi)- f(k)\big).
$$
The evolution of $\mcE^\ve(t)$ is given by
\begin{equation*}\begin{split}
&\partial_t \big<J,\mcE^\ve(t)\big>=
\gamma\frac \ve 2 \int_{\T} dk \; \big<|\hat{\psi}(k,t/\ve)|^2\big>_\ve CJ(k)\\
&\hspace{10pt}-\gamma\frac \ve 2 \int_{\T} dk \;\frac 1 2 [\big<(\hat{\psi}(k)\hat{\psi}(-k))(t/\ve)\big>_\ve+
\big<(\hat{\psi}(k)^*\hat{\psi}(-k)^*)(t/\ve)\big>_\ve](CJ)(k).
\end{split}\end{equation*}
Defining the distribution $Y^\ve(t)$ on $\T$ through
\begin{equation*}\begin{split}
\big<J,Y^\ve(t)\big>= & \frac \ve 2 \sum_{y,y'\in \Z}\big<\psi(y',
t/\ve)\psi(y,t/\ve)\big>_\ve
\int_\T dk\; e^{2\pi i k(y'-y)}J(k)\\
= &\frac \ve 2 \int_{\T} dk\;
\big<[\hat{\psi}(k)\hat{\psi}(-k)](t/\ve)\big>_\ve J(k),
\end{split}\end{equation*}
we can rewrite the evolution equation as
\begin{equation}\label{energyev}
\begin{split}
\partial_t \big<J,\mcE^\ve(t)\big>= \gamma\big<CJ, \mcE^\ve(t)\big> - \frac
\gamma 2 \big(\big<CJ,Y^\ve(t)\big>+\big<CJ,Y^\ve(t)^*\big>\big).
\end{split}
\end{equation}
This is not a closed equation for $\mcE^\ve(t)$. However we expect
that in the  limit $\ve\to 0$ the terms containing the distributions
$Y^\ve(t)$, $Y^\ve(t)^*$ disappear. In order to prove it, let us
consider the evolution of the distribution $Y^\ve(t)$ on the kinetic
time scale. Calculations are similar to the previous ones, but with
the difference that now $A [\hat{\psi}(k)\hat{\psi}(-k)]\neq 0$, and
precisely
$$
A [\hat{\psi}(k)\hat{\psi}(-k)]=-2i\omega (k)\hat{\psi}(k)\hat{\psi}(-k).
$$
We arrive at
\begin{equation}\label{ev:Y0}\begin{split}
\partial_t & \big<J,Y^\ve(t)\big>= -\frac{2i}{\ve}\big<\omega J,Y^\ve(t)\big>
+\frac \gamma 2 \big<\hat \beta J,Y^\ve(t)\big>\\
&+\frac \gamma 2 \big(\big<CJ,Y^\ve(t)\big>+\big<CJ,Y^\ve(t)^*\big>\big)
-\frac \gamma 2 \big<\hat \beta J,Y^\ve(t)^*\big>
- \gamma  \big<CJ, \mcE^\ve(t)\big>.
\end{split}\end{equation}
Observe that by integrating eq. (\ref{ev:Y0}) in time, we obtain
$$
\lim_{\ve\to 0}\Big|\int_0^tdt\;\big<\omega J,Y^\ve(t)\big> \Big|=0
$$
for every
bounded function $J$. In particular, since by item $(i)$ of lemma
\ref{lemma_omega}
$$
\sup_{k\in\T}\frac{R(k,k')}{\omega(k)}<\infty,
$$
we can choose a function $\omega^{-1}CJ$ with $J$ bounded  and  obtain
$$
\lim_{\ve\to 0}\Big|\int_0^tdt\;\big<C J,Y^\ve(t)\big> \Big|=0 .
$$
In the same way we have $\lim_{\ve\to 0}\big|\int_0^tdt\;\big<C
J,Y^\ve(t)^*\big> \big|=0$ and any limit distribution $\mcE(t)$ of
$\mcE^\ve(t)$ solves the equation
$$
\big<J,\mcE(t)\big>= \big<J, \mcE(0)\big> +
\gamma \int_0^t ds\; \big<CJ, \mcE(s)\big>
$$
for every bounded real valued function $J$.
$\qed$

\subsection{Proof of theorem \ref{th-inhbe}}
Now we will give the proof of (\ref{be}) for test functions
$J\in\mathcal{S}(\R\times\T, \C)$. The main difference to the
previous case is that  the Hamiltonian part of the generator
contributes to the evolution of $W^\ve(t)$, resulting in a ballistic
transport term. In order to control this term, we need to ensure
that there is no mass concentration at  $k=0$ for  every macroscopic
time $t\in[0,T]$ with $T>0$. This is stated in the following lemma.
%
\begin{lemma}\label{lemma_nomass}
Let assumption (b4) hold. Then for every $t\in[0,T]$
$$
\lim_{\rho\to 0}\,\overline{\lim_{\ve\to 0}} \; \frac \ve 2
\int_{|k|<\rho}dk \; \big<|\hat{\psi}(k,t/\ve)|^2\big>_\ve = 0.
$$
\end{lemma}
\emph{Proof.}
We use the evolution equation (\ref{energyev}) for $J_\rho(k) =
1_{[-\rho,\rho]}(k)$.
Since $|CJ_\rho(k)| \le c_1 (2\rho +
J_\rho(k))$
and $\big<|CJ_\rho(k)|,|Y^\ve(t)+ {Y^\ve(t)}^*|\big>\leq \big<|CJ_\rho(k)|,\mcE^\ve(t) \big>$,
we obtain the bound
\begin{equation*}
 \begin{split}
   \big<J_\rho, \mcE^\ve(t)\big> \le \big<J_\rho, \mcE^\ve(0)\big> +
   c_2\gamma \int_0^t ds\; \big< |CJ_\rho|, \mcE^\ve(s)\big>\\
   \le \big<J_\rho, \mcE^\ve(0)\big> +
   c_3\gamma\Big( 2\rho Kt + \int_0^t ds\;\big< J_\rho, \mcE^\ve(s)\big>  \Big),
 \end{split}
\end{equation*}
where $K$ is the bound on the total energy from condition (b3).
Then by Gronwall's inequality
\begin{equation*}
   \big<J_\rho, \mcE^\ve(t)\big> \le \left( 2\rho K +
     \big< J_\rho, \mcE^\ve(0)\big>  \right)e^{c_3\gamma t},
\end{equation*}
where, by assumption (b4), $\overline{\lim}_{\ve\to 0} \big< J_\rho,
\mcE^\ve(0)\big>\to 0$ for $\rho\to 0$. $\qed$


\subsubsection{Proof of  theorem \ref{th-inhbe}}
For every $J\in\mathcal{S}(\R\times\T,\C)$
the evolution of the distribution $W^\ve(t)$
on  the kinetic time-scale
is given by
\begin{equation*}
\begin{split}
\partial_t &\big< J ,W^\ve(t)\big> \\
&
=(\ve/2)\sum_{y,y'\in\Z}\partial_t\big<\psi(y',t/\ve)^*\psi(y,t/\ve)\big>_\ve
\int_\T dk\; e^{2\pi i k(y'-y)}J(\ve(y+y')/2,k)^*\\
& = (\ve/2)\sum_{y,y'\in\Z}\ve^{-1}\big<L[\psi(y')^*\psi(y)]\big>_\ve\int_\T dk\;
e^{2\pi i k(y'-y)}J(\ve(y+y')/2,k)^* ,
\end{split}\end{equation*}
where
\begin{equation*}\begin{split}
L[\psi(y')^*\psi(y)]=A[\psi(y')^*\psi(y)]+\ve\gamma S[\psi(y')^*\psi(y)]
\end{split}\end{equation*}
and $A$, $S$ are defined in (\ref{def:Agen}), (\ref{defq:Sgen1}),
respectively.
We start by computing the evolution determined by $A$,
the Hamiltonian part of the generator. Using the representation of
the Wigner distribution in
Fourier space we get
\begin{equation*}
\begin{split}
&\frac{\ve}2 \int_\R dp \int_\T dk\;  \ve^{-1}A\big[\big<\hat\psi(k-\ve p/2)^*
\hat\psi(k+\ve p/2)\big>_\ve ]
\widehat J(p,k)^* \\
&\hspace{10pt}=-i\frac{\ve}2 \int_\R dp \int_\T dk \big<\hat\psi(k-\ve p/2)^*
\hat\psi(k+\ve p/2)\big>_\ve \\
&\hspace{24pt} \times\ve^{-1}
[\omega(k+\ve p/2) - \omega(k-\ve p/2)]
\widehat J(p,k)^*.
\end{split}\end{equation*}
Now we prove that one can replace $\ve^{-1}[\omega(k+\ve p/2) -
\omega(k-\ve p/2)]$ with $\omega'(k)p$ in the last expression.
For every $0<\rho<1/2$
\begin{equation*}
\begin{split}
& \frac{\ve}2 \int_\R dp \int_\T dk \big<\hat\psi(k-\ve p/2)^*
\hat\psi(k+\ve p/2)\big>_\ve \\
&\qquad \times\left( \ve^{-1}
[\omega(k+\ve p/2) - \omega(k-\ve p/2)] -\omega'(k) p \right)
\widehat J(p,k)^* = I^\ve_>(\rho)+I^\ve_<(\rho),
\end{split}
\end{equation*}
where
\begin{equation*}
\begin{split}
I^\ve_>(\rho)=  & \frac{\ve}2 \int_\R dp \int_{|k|>\rho} dk
\;\big<\hat\psi(k-\ve p/2)^* 
\hat\psi(k+\ve p/2)\big>_\ve \\
&\qquad \times\left( \ve^{-1}
[\omega(k+\ve p/2) - \omega(k-\ve p/2)] -\omega'(k) p \right)
\widehat J(p,k)^*,\\
I^\ve_<(\rho)=  & \frac{\ve}2 \int_\R dp \int_{|k|<\rho} dk
\;\big<\hat\psi(k-\ve p/2)^* 
\hat\psi(k+\ve p/2)\big>_\ve \\
&\qquad \times\left( \ve^{-1}
[\omega(k+\ve p/2) - \omega(k-\ve p/2)] -\omega'(k) p \right)
\widehat J(p,k)^*.
\end{split}
\end{equation*}
Using Schwarz inequality and points (i), (ii) of lemma
\ref{lemma_omega} in the appendix
\begin{equation*}\begin{split}
\big| I^\ve_<(\rho)\big|&\leq \int_\R dp\;
(C|p|+\|\nabla\omega\|_\infty)
\sup_{k\in\T}|\widehat J(p,k)|\Big(
\frac{\ve}{2}\int_{|k|\leq\rho}dk\;\big<|\hat{\psi}(k)|^2\big>_\ve
\Big)\\
&\leq C_0 \frac{\ve}{2}\int_{|k|\leq\rho}dk\;\big<|\hat{\psi}(k)|^2\big>_\ve
\end{split}\end{equation*}
and, by lemma \ref{lemma_nomass},
\begin{equation*}
\lim_{\rho\to 0}\overline{\lim_{\ve\to 0}}\;\big| I^\ve_<(\rho)\big|
= 0 .
\end{equation*}
To compute $I^\ve_>(\rho)$ we split it into two parts,
\begin{equation*}\begin{split}
I^\ve_>(\rho)=  & \frac{\ve}2 \int_{\R,\ve|p|\geq\rho} dp
\int_{|k|>\rho} dk \;\big<\hat\psi(k-\ve p/2)^* 
\hat\psi(k+\ve p/2)\big>_\ve \\
&\qquad \times\left( \ve^{-1}
[\omega(k+\ve p/2) - \omega(k-\ve p/2)] -\omega'(k) p \right)
\widehat J(p,k)^*\\
& + \frac{\ve}2 \int_{\R,\ve|p|<\rho} dp \int_{|k|>\rho} dk\;
\big<\hat\psi(k-\ve p/2)^* 
\hat\psi(k+\ve p/2)\big>_\ve \\
&\qquad \times\left( \ve^{-1}
[\omega(k+\ve p/2) - \omega(k-\ve p/2)] -\omega'(k) p \right)
\widehat J(p,k)^*.
\end{split}\end{equation*}
Again we use Schwarz inequality and points (i), (ii) of lemma
\ref{lemma_omega} to show that the first term on the RHS is
negligible, since for all $\rho >0$ it is bounded by
\begin{equation*}
K\int_{|p|\geq\rho/\ve}dp\;(c|p|+\|\nabla\omega\|_\infty)\sup_{k\in\T}|\widehat
J(p,k)|,
\end{equation*}
which tends to 0 as $\ve \to 0$.

For the second term on the RHS we use the point (iii) of lemma
\ref{lemma_omega}, since $|k|>\rho$, $\ve|p|<\rho$ implies
$|k|>\ve|p|$, and for all $\rho >0$ we get
\begin{equation*}\begin{split}
&\frac{\ve}2 \int_{|p|<\rho/\ve} dp \int_{|k|>\rho} dk \;\big<\hat\psi(k-\ve p/2)^*
\hat\psi(k+\ve p/2)\big>_\ve \\
&\qquad\times\left( \ve^{-1}
[\omega(k+\ve p/2) - \omega(k-\ve p/2)] -\omega'(k) p \right)
\widehat J(p,k)^*\\
& \leq K \int_{\R}dp\;
\ve \frac{C_4}{\rho}|p|^2\sup_{k\in\T}|\widehat{J}
(p,k)|,
\end{split}\end{equation*}
which tends to 0 as $\ve \to 0$.
Then we have
\begin{equation}\label{hampart}
\begin{split}
& \frac{\ve}2 \int_\R dp \int_\T dk\; \ve^{-1}A\big[\big<\hat\psi(k-\ve p/2)^*
\hat\psi(k+\ve p/2)\big>_\ve ]
\widehat J(p,k)^* \\
&\hspace{10pt}= \frac{\ve}2 \int_\R dp \int_\T dk\; \big<\hat\psi(k-\ve p/2)^*
\hat\psi(k+\ve p/2)\big>_\ve (-i\;p)
\omega'(k)\widehat J(p,k)^* + \mathcal{O}(\ve) \\
&\hspace{10pt} =\frac{1}{2\pi}\Big< \nabla\omega \nabla_r J, W^\ve(t)\Big>_\ve
+\mathcal{O}(\ve).
\end{split}
\end{equation}

We have to compute
\begin{equation*}
\begin{split}
&(\ve/2)\sum_{y,y'\in\Z}\big<S[\psi(y')^*\psi(y)]\big>_\ve\int_\T dk\;
e^{2\pi i k(y'-y)}J(\ve(y+y')/2,k)^*\\
&\hspace{10pt}= \frac{\ve}{2}\sum_{y,y'\in\Z}\big<S[\psi(y')^*\psi(y)]\big>_\ve
\tilde J\big(\ve (y+y')/2, y - y'\big)^*.
\end{split}
\end{equation*}
Since $S$ is a second order operator, we have
\begin{equation*}
\begin{split}
S[\psi(y')^*\psi(y)]=\psi(y')^*S\psi(y)+[S\psi(y')^*]\psi(y)
+\frac{1}{3}\sum_{z\in\Z}[Y_z\psi(y')^*][Y_z\psi(y)],
\end{split}
\end{equation*}
where by direct computation
\begin{equation}\label{S2}
S\psi(y) = \frac{i}{\sqrt 2} S p_y = \frac 12 \beta * (\psi^* -
\psi) (y) 
\end{equation}
and
\begin{equation*}
\begin{split}
\psi(y')^*S\psi(y)+[S\psi(y')^*]\psi(y)=\frac{i}{\sqrt{2}}
\big(\psi(y')^*(\beta *p)_y-\psi(y)(\beta * p)_{y'}\big).
\end{split}
\end{equation*}
Integrating by parts and using the symmetry of $\beta$,
we can rewrite
\begin{equation*}
\begin{split}
&\gamma (\ve/2)\sum_{y,y'\in\Z}
\big[\big<\psi(y')^*S\psi(y)\big>_\ve+\big<(S\psi(y')^*)\psi(y)\big>_\ve\big]
\tilde
J\big(\ve (y+y')/2, y - y'\big)^*\\
&\hspace{10pt}= \gamma (\ve/2)\sum_{y,y'\in\Z}\frac{i}{\sqrt{2}}\left< \psi(y')^*
p_y\right>_\ve \sum_{z\in\Z}\beta(z)
\tilde J\big(\ve (y+y'+z)/2, y - y'+ z\big)^*\\
&\hspace{10pt}-\gamma (\ve/2)\sum_{y,y'\in\Z}\frac{i}{\sqrt{2}}\left< \psi(y)
p_{y'}\right>_\ve \sum_{z\in\Z}\beta(z) \tilde J\big(\ve (y+y'-z)/2,
y - y'+ z\big)^*.
\end{split}
\end{equation*}
Using the energy bound and the properties of the test functions $J$,
one can write the first term on right hand side as
\begin{equation*}
\begin{split}
&\gamma (\ve/2)\sum_{y,y'\in\Z}\frac{i}{\sqrt{2}}\left<
\psi(y')^*p_y\right>_\ve \sum_{z\in\Z}\beta(z)\tilde
J\big(\ve (y+y')/2, y - y'+ z\big)^* +\mathcal{O}(\ve) \\
&\hspace{10pt}= \gamma (\ve/2)\sum_{y,y'\in\Z}\frac{i}{\sqrt{2}}\left<
\psi(y')^*p_y\right>_\ve \int_\T dk e^{2\pi i k (y'-
y)}\hat{\beta}(k) J\big(\ve\ (y+y')/2, k\big)^* \\
&\hspace{24pt}+\mathcal{O}(\ve)
\end{split}
\end{equation*}
and the same can be done for the other term.
Finally we obtain
\begin{equation}\label{eq:S2}
\begin{split}
&(\ve/2)\sum_{y,y'\in\Z}
\big[\big<\psi(y')^*S\psi(y)\big>_\ve+\big<(S\psi(y')^*)\psi(y)\big>_\ve\big]
\tilde
J\big(\ve (y+y')/2, y - y'\big)^*\\
&\hspace{10pt}=(\ve/2)\sum_{y,y'\in\Z}\frac{i}{\sqrt{2}}
\big<\psi(y')^*p_y-\psi(y)p_{y'}\big>_\ve
\int_\T dk \;e^{2\pi i k (y'- y)}\hat{\beta}(k) J\big(\ve (y+y')/2,
k\big)^* \\ &\hspace{280pt}+\mathcal{O}(\ve).
\end{split}
\end{equation}

About the other term in (\ref{S2}) first observe that $
\sum_{z\in\Z} [Y_z\psi(y')^*][Y_z\psi(y)]$ is just a finite sum for
any $y,y'$ fixed and, computing it explicitly and identifying terms
that differ by translations, see details of the computation in
appendix \ref{proof:SW}, one obtains
\begin{equation}\label{eq:Ypsi}
\begin{split}
&\frac{1}{3}(\ve/2) \sum_{y,y',z\in\Z} \big< [Y_z\psi(y')^*][Y_z\psi(y)]\big> \tilde
J(\ve (y'+y)/2, y - y')^*\\
&\hspace{10pt}= (\ve/2) \sum_{y\in\Z} \sum_{z,u=-2}^2 \alpha(z,u) \big<p_y p_{y+z}\big>_\ve
\tilde J( \ve y , u)^* + \mathcal{O}(\ve),
\end{split}\end{equation}
where $\alpha(z,u)=\alpha(-z,u)=\alpha(z,-u)=\alpha(u,z)$ and is given by
(\ref{eq:alpha}).
We can rewrite it as
\begin{equation*}
\begin{split}
&(\ve/2) \sum_{y\in\Z} \sum_{z,u=-2}^2 \alpha(z,u) \big<p_y p_{y+z}\big>_\ve
\tilde J( \ve( y+ z/2) , u)^* + \mathcal{O}(\ve) \\
&\hspace{10pt}=(\ve/2) \sum_{y\in\Z} \sum_{z,u\in\Z} \alpha(z,u) \big<p_y
p_{y+z}\big>_\ve \tilde J( \ve( y+z/2) , u)^* + \mathcal{O}(\ve),
\end{split}\end{equation*}
where we put $\alpha(z,u)=0$ if $|z|>2$ or $|u|>2$ and by changing
variables we obtain
\begin{equation}\label{45}\begin{split}
(\ve/2) \sum_{y,y'\in\Z} \big<p_y p_{y'}\big>_\ve \sum_{u\in\Z}
\alpha(y' - y,u) \tilde J\left( \ve (y'+y)/2 , u\right)^* +
\mathcal{O}(\ve).
\end{split}
\end{equation}
Defining
\begin{equation*}
R(k,k') = \sum_{z\in \Z} \sum_{u\in\Z}e^{- 2\pi i k z}e^{-2\pi i k'
u} \alpha(z, u),
\end{equation*}
we can rewrite (\ref{45}) as
\begin{equation*}
\begin{split}
(\ve/2) \sum_{y,y'\in\Z} \big<p_y p_{y'}\big>_\ve \int dk\; e^{2\pi
i k(y'-y)} \int dk'\; R(k,k') J(\ve (y'+y)/2 , k')^* +
\mathcal{O}(\ve),
\end{split}
\end{equation*}
where direct computation gives
\begin{equation*}
\begin{split}
R(k,k') = \frac 23\big( 3 - 2 \cos(2\pi k) - \cos(4\pi k) - 2
\cos(2\pi k')
+ 2\cos(2\pi (k'+ 2k))\\ - \cos(4\pi k') + 2
\cos(2\pi(2k'+ k)) - \cos(2\pi(2k'+ 2k)) \big)\\
=\frac{4}{3}\big(2\sin^2(2\pi k)\sin^2(\pi k')+
2\sin^2(2\pi k')\sin^2(\pi k)
-\sin^2(2\pi k)\sin^2(2\pi k')\big),
\end{split}
\end{equation*}
which is the kernel defined in (\ref{kernel1}).
Using the relation $\int_\T dk' R(k,k')=-\hat{\beta}(k)$ and
\begin{equation*}\begin{split}
&\big<p_yp_{y'}\big>_\ve\\
&\hspace{10pt}= \frac{1}{2}[\big<\psi(y')^*\psi(y)\big>_\ve
+\big<\psi(y')\psi(y)^*\big>_\ve]
-\frac{1}{2}[\big<\psi(y')\psi(y)\big>_\ve
+\big<\psi(y')^*\psi(y)^*\big>_\ve]\\&
\frac{i}{\sqrt{2}}\big<\psi(y')^*p_y-\psi(y)p_{y'} \big>_\ve\\
&\hspace{10pt}=
\big<\psi(y')^*\psi(y)\big>_\ve
-\frac{1}{2}[\big<\psi(y')\psi(y)\big>_\ve
+\big<\psi(y')^*\psi(y)^*\big>_\ve],
\end{split}\end{equation*}
we obtain
\begin{equation}
\begin{split}
\gamma \frac{\ve}{2} \sum_{y,y'\in\Z} \left< S\left(\psi(y')^*
\psi(y)\right)\right>_\ve J(\ve (y'+y)/2, y' - y)\hspace{1cm}
\\= \gamma\left<CJ, W^\ve\right> - \frac {\gamma}{2} (\left<CJ,
Y^\ve\right> + \left<CJ, Y^{\ve*}\right>)+\mathcal{O}(\ve),
\end{split}
\end{equation}
where the \emph{collision operator} $C$ is defined in (\ref{C1}) and
the distributions $Y^\ve(t)$, $Y^\ve(t)^*$ are  defined as
\begin{equation}\label{def:Y,Y*}\begin{split}
\big<J,Y^\ve(t)\big>= &
(\ve/2)\sum_{y,y'\in\Z}\left<
\psi(y') \psi(y)\right>_\ve
\tilde J(\ve (y+y')/2, y - y')^* ,\\
\big<J,Y^\ve(t)^*\big>= & (\ve/2)\sum_{y,y'\in\Z}\left< \psi(y')^*
\psi(y)^*\right>_\ve \tilde J(\ve (y+y')/2, y - y')^*
\end{split}\end{equation}
for every $J\in\mathcal{S}(\R\times\T,\C)$.

The evolution of $W^\ve(t)$ is not a closed equation,
\begin{equation}\label{evol:W}
\begin{split}
\partial_t\big<J,W^\ve(t) \big>= & \big<(\nabla\omega \nabla_x
J),W^\ve(t)\big> +\gamma
\big<(C J),W^\ve(t)\big>\\ &
-\frac{\gamma}{2}\big<(CJ),Y^\ve(t)\big>-
\frac{\gamma}{2}\big<(CJ),Y^\ve (t)^*\big> + \mathcal{O}(\ve).
\end{split}
\end{equation}
However we expect that in the kinetic limit $\ve\to 0 $ the terms
containing the distributions $Y^\ve(t), Y^\ve(t)^*$ to disappear. To
prove this, we consider the evolution of $Y^\ve(t)$. Again by Ito's
formula
\begin{equation*}
\begin{split}
\partial_t \big< J ,Y^\ve(t)\big> &=
(\ve/2)\sum_{y,y'\in\Z}\partial_t\big<\psi(y')\psi(y)\big>_\ve
\tilde J(\ve(y+y')/2,y-y')^*\\
&= (\ve/2)\sum_{y,y'\in\Z}\ve^{-1}\big<L[\psi(y')\psi(y)]\big>_\ve\tilde J(\ve(y+y')/2,y-y')^*\\
&=
(\ve/2)\sum_{y,y'\in\Z}\ve^{-1}\big<A[\psi(y')\psi(y)]\big>_\ve\tilde
J(\ve(y+y')/2,y-y')^* \\
&+\gamma (\ve/2)\sum_{y,y'\in\Z}\big<S[\psi(y')\psi(y)]\big>_\ve\tilde
J(\ve(y+y')/2,y-y')^* .
\end{split}\end{equation*}
For the stochastic part of the generator, by a similar
computation as above, we obtain
\begin{equation}\label{eq:SY}
\begin{split}
&\gamma (\ve/2)\sum_{y,y'\in\Z}\big<S[\psi(y')\psi(y)]\big>_\ve \tilde
J(\ve(y+y')/2,y'-y) \\
&\hspace{10pt}=\frac{\gamma}{2}\big<(CJ),Y^\ve\big> +
\frac{\gamma}{2}\big<(CJ),{Y^\ve}^*\big>+\frac{\gamma}{2}\big<(\hat{\beta}J),Y^\ve\big>
-\frac{\gamma}{2}\big<(\hat{\beta}J),{Y^\ve}^*\big>\\
&\hspace{24pt}-\frac{\gamma}{2}\big[\big<(CJ),W^\ve\big>+\big<(CJ),{W^\ve}^*\big>\big]+\mathcal{O}(\ve).
\end{split}
\end{equation}


To compute the Hamiltonian contribution to the evolution of
$Y^\ve(t)$, we use the representation of $Y^\ve(t)$ in the Fourier space
and get
\begin{equation*}
\begin{split}
&(\ve/2)\int_{\R}dp\int_\T dk\;\ve^{-1}
A[\big<\hat{\psi}(k-\ve p/2) \hat{\psi}(k+\ve p/2)\big>_\ve]\\
&\hspace{10pt}=-2i \ve^{-1}(\ve/2)\int_{\R}dp\int_\T dk\;\ve^{-1}
\big<\hat{\psi}(k-\ve p/2) \hat{\psi}(k+\ve p/2)\big>_\ve\\
&\hspace{60pt}(\omega(k+\ve
p/2) + \omega(k- \ve p/2))\widehat{J}(p,k)^*,
\end{split}\end{equation*}
where with similar arguments as  above  one  can replace
$\omega(k+\ve p/2) + \omega(k- \ve p/2)$ by $2\omega(k)$, with an
error of order $\ve^2$. Then we arrive at the following equation for
the evolution of $Y^\ve(t)$,
\begin{equation}\label{evol:Y}
\begin{split}
&\partial_t\big<J,Y^\ve(t) \big> =
-\frac{2i}{\ve}\big<(\omega J),Y^\ve(t)\big> +
\frac{\gamma}{2}\big<(CJ),Y^\ve\big>+\frac{\gamma}{2}\big<(\hat{\beta}J),Y^\ve\big>\\& +
\frac{\gamma}{2}\big<(CJ),{Y^\ve}^*\big>
-\frac{\gamma}{2}\big<(\hat{\beta}J),{Y^\ve}^*\big>
-\frac{\gamma}{2}\big[\big<(CJ),W^\ve\big>+\big<(CJ),{W^\ve}^*\big>\big]+\mathcal{O}(\ve).
\end{split}
\end{equation}
After time integration, we obtain, for any $J \in\mathcal{S}(\mathbb R\times \mathbb T)$,
\begin{equation*}
\lim_{\ve \to 0} \Big| \int_0^t dt\; \big<(\omega J),Y^\ve(t)\big>\Big| = 0.
\end{equation*}
Observe that $R(k,k')/\omega(k)\in C^\infty(\T/\{0\})$ and,
using item (i) of lemma \ref{lemma_omega} in the appendix,
$$\sup_{k\in\T} \frac{R(k,k')}{\omega(k)} < \infty .$$
Then equation (\ref{evol:Y}) holds for any function
\begin{equation*}
\omega(k)^{-1} C J (y,k) = \int_\T dq\;\frac{R(k,q)}{\omega(k)}
\left( J(y, q) - J(y,k)\right) 
\end{equation*}
with $J\in\mathcal{S}(\R\times\T)$
and consequently
\begin{equation*}
\lim_{\ve \to 0} \Big| \int_0^t dt\;\big<( C J),Y^\ve(t)\big>\Big| = 0.
\end{equation*}


\section{Extension to  dimensions $d\ge 2$}
\label{sec:extens-d-dimens}

We consider a particular generalisation of our model to $d$
dimensions, $d\geq 2$. The perfect lattice is $\Z^d$. Deviations
from the equilibrium position $\by\in\Z^d$ is $\bq_{\by}\in\R^d$ and
$\bp_{\by}$ denotes the corresponding momentum. Thus the phase space
is $(\R^d\times\R^d)^{\Z^d}$. The Hamiltonian of the system is given
by
\begin{equation}
\label{def:Hamiltd}
H(\bp,\bq) = \frac 12 \sum_{\by\in\Z^d} {\bp_\by^2} +
\frac 12 \sum_{\by,\by'\in\Z^d} \alpha(\by - {\by'})
\bq_\by \cdot \bq_{\by'}.
\end{equation}
For simplicity the couplings $\alpha$ are taken to be scalar. In
general,  $\alpha$ would be a $d\times d$ matrix.
We denote
\begin{equation}
\hat{v}(\bk)=\sum_{\bz\in\Z^d}e^{-2\pi i\bk\cdot\bz} v(\bz), \qquad
\tilde{f}(\bz)=\int_{\T^d}d\bk \; e^{2\pi i\bk\cdot\bz} f(\bk).
\end{equation}
We assume 
$\alpha(\cdot)$ to satisfy the following properties:
\begin{assumption}\label{alpha_assumpd}\hspace{4cm}
\begin{itemize}
\item\emph{(a1)} $\alpha(\by)\neq 0$ for some $\by\neq 0$.
\item\emph{(a2)} $\alpha(\by)=\alpha(-\by)$ for all $\by\in\Z^d$.
\item\emph{(a3)} There are constants $C_1,C_2>0$ such that for all $\by$
$$
|\alpha(\by)|\leq C_1e^{-C_2|\by|}.
$$
\item\emph{(a4)}
  \begin{itemize}
  \item \emph{(pinning):} $\hat{\alpha}>0$ on $\T^d$,
    \item \emph{(no pinning):}
    $\hat{\alpha}(\bk)>0$ for all $\bk\neq 0$, $\hat{\alpha}(0)=0$,
    $\text{Hess}(\hat{\alpha}(0))$ is invertible.
  \end{itemize}

\end{itemize}
\end{assumption}

The dynamics is determined by the generator $L=A+\ve\gamma S$
with
\begin{equation}
\label{def:Agend}
\begin{split}
A = \sum_{\by\in\Z^d} \bp_{\by} \cdot \partial_{\bq_{\by}} -
\sum_{\by,\by'\in\Z^d} \alpha(\by - \by')
\bq_{\by'} \cdot \partial_{\bp_{\by}} .
\end{split}
\end{equation}
$S$ is defined  through the vector fields
\begin{equation*}
\label{eq:Xfield}
X^{i,j}_{\bx, \bz} = (p^j_\bz-p^j_\by) (\partial_{p^i_\bz}
- \partial_{p^i_\by}) -(p^i_\bz-p^i_\by) (\partial_{p^j_\bz}
- \partial_{p^j_\by}) ,
\end{equation*}
according to
\begin{equation}
\begin{split}
\label{def:Sgen}
S  = \displaystyle\frac 1{2(d-1)} \sum_{\by\in\Z^d}\sum_{i,j,k=1}^d
\big( X^{i,j}_{\by, \by+\be_k}\big)^2 
= \displaystyle\frac 1{4(d-1)} \sum_{\by, \bz \in \Z^d, \atop \|\by -
\bz\|=1}\sum_{i,j=1}^{d}
\left( X^{i,j}_{\bx, \bz}\right)^2,
\end{split}\end{equation}
where ${\be}_1,\ldots,{\be}_d$ is the canonical basis of ${\mathbb
Z}^d$. As in the one-dimensional case
\begin{equation*}
S \; \sum_{\by\in\Z^d} \bp_\by = 0\ ,\hspace{0.4cm} S H = 0 .
\end{equation*}
Note that now it suffices to couple nearest neighbors.

The evolution of $\{\bp(t),\bq(t)\}$ is given by the following
stochastic differential equations
\begin{equation}
\label{eq:sde}
\begin{split}
d\bq_\by &= \bp_\by\; dt,\\
d\bp_\by &= -(\alpha * \bq)_{\by}\; dt +
2 \ve\gamma \Delta \bp_\by \; dt \\
& \qquad  + \frac {\sqrt{\ve\gamma}}{2\sqrt{ d-1}}
\sum_{\bz \in\Z^d, \atop \|\bz - \by\|=1} \sum_{i,j=1}^d
\left(X^{i,j}_{\by, \bz} \bp_{\by} \right) \; dw^{i,j}_{\by, \bz}(t)
\end{split}
\end{equation}
for all $\by\in\Z^d$.
Here
$\{w^{i,j}_{\bz, \by} = w^{i,j}_{\by, \bz};\; \bz, \by \in \Z^d;\;
i,j= 1,\dots, d;\; \|\by - \bz\| = 1\}$ are independent standard
Wiener processes.

As before, we define the complex valued vector field  $\bpsi:\Z^d\to\C^d$ by
\begin{equation}\label{def:psid}\begin{split}
\bpsi(\by,t)=\frac{1}{\sqrt{2}}\big((\tilde{\omega}\ast \bq)_\by(t)+i
\bp_\by(t)\big)
\end{split}\end{equation}
with the inverse relation
\begin{equation}
\bp_\by(t)=\frac{i}{\sqrt 2}(\bpsi^*-\bpsi)(\by,t).
\end{equation}
Observe that $|\bpsi(\by)|^2 = e_\by$, the local energy ay $\by$.
For every $t\geq 0$, the evolution of $\bpsi$ is given by the
stochastic differential equation,
\begin{equation}\begin{split}
d\bpsi(\by,t)= & -i(\tom\ast \bpsi)(\by,t) dt
+\frac{1}{2}\ve\gamma\beta*(\bpsi -
\bpsi^*)(\by,t) dt\\
& +\frac {\sqrt{\ve\gamma}}{4\sqrt{ d-1}}
\sum_{\by'\in\Z^d, \atop \|\by' - \by\|=1} \sum_{i,j=1}^d \big(X^{i,j}_{\by, \by'}
(\bpsi-\bpsi^*)(\by,t) \big) \; dw^{i,j}_{\by, \by'}(t) ,
\end{split}\end{equation}
where $\beta$ is  determined through $(\beta * f)(\bz)= \Delta f(\bz)$.

Given a function $J$ on $\R^d\times\T^d$, we define
\begin{equation}
\label{eq:5d}
\tilde J(\bx, \bz) = \int_{\T^d}d\bk\; e^{i2\pi\bk\cdot\bz} J(\bx,\bk) 
\end{equation}
on $\R^d\times\Z^d$.
We also define
\begin{equation}
\widehat J(\bp,\bk) = \int_{\R^d}
d\bx\; e^{-i2\pi\bp\cdot\bx} J(\bx,\bk)  .
\end{equation}
We choose a class of test-functions $J$ on
$\R^d\times\T^d$ such that $J\in
\mathcal{S}(\R^d\times\T^d,\mathbb{M}_d)$, where $\mathbb{M}_d$ is the
space of complex $d\times d$ matrices.

Fix $\ve>0$. We introduce the complex valued correlation matrices
\begin{equation}\label{corr_matrix}\begin{split}
\big<\bpsi(\by')^*\otimes\bpsi(\by)\big>_\ve, \qquad
\big<\bpsi(\by')\otimes\bpsi(\by)\big>_\ve,
\end{split}\end{equation}
where $\big<\cdot\big>_\ve$ denotes the expectation value with respect
to a probability measure  on phase  space which  satisfies the
following properties:
\begin{enumerate}[(c1)]
\item $\big< \bpsi(\by)\big>_\ve=0, \hspace{0.4cm}\forall\by\in\Z^d,$
\item $ \big< \bpsi(\by')\otimes\bpsi(\by)\big>_\ve=0, \hspace{0.4cm}\forall\by, \by'\in\Z^d,$
\item $ \big< \|\bpsi\|^2\big>_\ve=\big<\sum_{\bz\in\Z^d}|\bpsi(\bz)|^2\big>_\ve\leq K\ve^{-d}.$
\end{enumerate}
Observe that, since $\big< \|\bpsi\|^2\big>_\ve=\big< H \big>_\ve$, we are considering states with
an energy of order $\ve^{-d}$.
We define the matrix-valued Wigner distribution $W^\ve$ as
\begin{equation}
\begin{split}
&\big<J,W^\ve \big> \\&\hspace{6pt} =(\ve/2)^d\sum_{\by, \by'\in\Z^d}\sum_{i,j=1}^d
\big<\psi_j(\by')^*\psi_i(\by)\big>_\ve \int_{\T^d} d\bk \; e^{i2\pi
\bk\cdot(\by'-\by)} J_{j,i}(\ve(\by'+\by)/2,k)^*\\
&\hspace{6pt} = (\ve/2)^d\sum_{\by, \by'\in\Z^d}\sum_{i,j=1}^d
\big<\psi_j(\by')^*\psi_i(\by)\big>_\ve \tilde
J_{j,i}(\ve(\by'+\by)/2,\by-\by')^*
\end{split}
\end{equation}
with
$J\in \mathcal S (\R^d\times\T^d, \mathbb{M}_d)$.
The evolution of the diagonal terms of the
distribution $W^\ve$ on  time scale $\ve^{-1}t$ is determined through
\begin{equation}
\begin{split}
&\big<J,W^\ve(t) \big>
=(\ve/2)^d\sum_{i=1}^d \sum_{\by, \by'\in\Z^d}\big<\psi_i(\by',
t/\ve)^*\psi_i(\by, t/\ve)\big>_\ve\\
&\hspace{100pt}\times\int_{\T^d} d\bk \;
e^{i2\pi\bk\cdot(\by'-\by)} J_{i}(\ve(\by'+\by)/2,k)^*
\end{split}
\end{equation}
for $J_i \in \mathcal{S}(\R^d\times\T^d)$, $i=1\dots,d$.

Observe that since the dynamics preserves the total energy, the
condition $\ve^d \left<\|\psi \|\right> \le K$ holds at any time 
and, by proposition \ref{prop:bound}, the Wigner distribution is
well defined at any time. On this time scale the diagonal terms of
the distribution $W^\ve(t)$ converge in a weak sense to a (vector
valued) measure $\mu = \{\mu_i(t), i=1,..,d\}$ on $\R^d\times\T^d$
which satisfies the following Boltzmann equation. For any vector
valued function $J\in\mathcal S(\R^d\times\T^d,\C^d)$,
\begin{equation}
\label{Bed}
\begin{split}
\big< J, \mu(t) \big> -\big< J, \mu(0) \big> = \frac{1}{2\pi}\int_0^t
ds \;\left( \big<
\nabla\omega\cdot \nabla_{\bx} J, \mu(s) \big>
+ \gamma \big< C J, \mu(s) \big> \right) ,
\end{split}
\end{equation}
where $\big< J, \mu(t) \big>$ denotes the scalar product
$\sum_{i=1}^d\int_{\R^d\times\T^d} J_i(\bx,\bk)^*\mu_i(d\bx,d\bk)$.
The collision operator is given by
\begin{equation}\label{Coll}
\begin{split}
(CJ)_i (\bx, \bk) = \frac{1}{d-1}
\sum_{1\leq j\leq d, \atop j\neq i}\int_{\T^d}d\bk'\;R(\bk,\bk')
\big(J_j(\bx, \bk')- J_i(\bx, \bk)\big),
\end{split}
\end{equation}
where the kernel $R:\T^d\times\T^d\to\R$ has the following expression,
\begin{equation}\label{kernel}
R(\bk,\bk')=16\sum_{\ell=1}^d\sin^2(\pi k_\ell)\sin^2(\pi
{k'}_\ell).
\end{equation}
As in the one-dimensional case, in order to prove the inhomogeneous
Boltzmann equation (\ref{Bed}), we need an additional condition on
the initial distribution in the unpinned case $(\hat{\alpha}(0)=0)$,
which
ensures that there is no initial concentration of energy at $\bk=0$:\smallskip\\
(c4) In the case of no pinning we require
$$
\lim_{R\to 0}\,\overline{\lim_{\ve\to
    0}}\;(\ve/2)^d\int_{|\bk|<R}d\bk\;\big<|\hat{\bpsi}(\bk)|^2\big>_\ve =0.
$$
Now we state the precise theorem. The proof is analogous to the
one-dimensional case.
\begin{theorem}
Let \emph{Assumptions (c1-c4)} hold and assume that $W^\ve(0)$
converges to a positive vector valued measure $\mu_0$. Then, for all
$t\in[0,T]$, $W^\ve(t)$ converges to a positive (vector-valued)
measure $\mu_0(t)$ which is the unique solution of the Boltzmann
equation
\begin{equation}
\begin{split}
\partial_t \big< J, \mu(t) \big> =
\frac{1}{2\pi}\big< \nabla\omega\cdot \nabla_{\bx} J, \mu(t)
\big> + \gamma \big< C J, \mu(t) \big>
\end{split}
\end{equation}
with initial condition $\mu_0(t)$.
\end{theorem}

As in the one-dimensional case, the Boltzmann equation has a
probabilistic interpretation as the forward equation of a Markov
process. We consider the Markov process
$$\big(\bX(t),\bK(t),i(t)\big).$$
By (\ref{Coll}), the jump rate from $(i,\bk)$ to $(j, d\bk')$ is
given by
$$\nu_{\bk,i}(j,d\bk')=\frac{1}{d-1}(1-\delta_{i,j})R(\bk,\bk')d\bk',\hspace{0.4cm}\forall
i,j=1,...,d .$$Transitions between states with the same index
$i$ are forbidden. The total collision rate is
$$
\phi_i(\bk)=\sum_{j=1}^d\int_{\T^d}\nu_{\bk,i}(j,d\bk')
=\int_{\T^d}d\bk'\; R(\bk,\bk'),
$$
$i=1,...,d$, which does not  depend on $i$. Explicitly
\begin{equation}\label{phi}
\phi_i(\bk)=\phi(\bk)=8\sum_{\ell=1}^d\sin^2(\pi k_\ell).
\end{equation}
Given a state $(\bk,i)$ at $t=0$, it jumps at time $\tau$ to the
state $(d\bk',j)$ with a probability
$\nu_{\bk,i}(d\bk',j)/\phi(\bk)$, where $\tau$ is an exponentially
distributed random variable of mean $\phi(\bk)^{-1}$. As before the
position process $\bX(t)$ is defined through
\begin{equation}
\bX(t)=\bX(0)+ \frac{1}{2\pi}\int_0^t ds\;\nabla\omega(\bK(s)) .
\end{equation}


\section{Homogeneous case: correlations, energy current and
 conductivity}
\label{sec:homog-case:-corr}

\subsection{Translation invariant measures}
\label{sec:transl-invar-meas}

We consider a situation where the initial measure on phase space is
invariant under space translations. For simplicity we work in the
one-dimensional setting. Using the methods from section
\ref{sec:extens-d-dimens}, the generalization to $d\ge 2$  is
straightforward.

Since energy will be now a.s. infinite, the results of section
\ref{wignerd} do not apply. Let us denote with $\big<\cdot\big>$ the
expectation value with respect to this initial translation
invariant measure, and assume that it has the following properties:
\begin{enumerate}[(d1)]
\item
$\big<\psi(y)\big>= 0, \hspace{1cm} \forall y\in\Z$,
\item
$\big<\psi(y)\psi(0)\big>=0,\hspace{1cm}\forall y\in\Z$,
\item
$\sum_{z\in\Z} \big|\big<\psi(0)^* \psi(z)\big>\big|<\infty$.
\end{enumerate}

The Wigner distribution is still well defined for every function
$J\in\mathcal{S}(\R\times\T)$. Using translation invariance
\begin{equation*}\begin{split}
\big<J,W^\ve \big>&=\frac\ve 2 \sum_{y, y'\in\Z}
\big<\psi(y')^*\psi(y)\big> \int_{\T} dk e^{2\pi i
k(y'- y)} J (\ve(y'+y)/2,k)^*\\
&= \frac\ve 2 \sum_{y, y'\in\Z} \big<\psi(0)^*\psi(y-y')\big>
\int_{\T} dk e^{2\pi i
k(y'- y)} J (\ve(y'+y)/2,k)^*\\
&= \sum_{z\in\Z}\big<\psi(0)^*\psi(z)\big>\int_\T dk e^{-2\pi i k z}
\big[\frac\ve 2 \sum_{y\in\Z} J (\ve(2y+z)/2,k)^*\big] ,
\end{split}\end{equation*}
which is finite by condition (d3) and by the fast decay of $J$.
In Fourier space the previous expression becomes
\begin{equation*}
\begin{split}
&\big<J,W^\ve \big>\\
&\hspace{10pt}=\frac{\ve}{2}\sum_{y, y'\in\Z} \int_\R
dp\;\int_\T dk\; e^{2\pi i(k-\ve p/2)y'}\big<\psi(y')^*\psi(y)\big>
e^{-2\pi i(k+\ve p/2)y}\widehat{J}(p,k)^*
\end{split}\end{equation*}
and using  translation invariance
\begin{equation*}\begin{split}
\big<J,W^\ve \big>=\frac{1}{2}\int_\R dp\;\int_\T dk\;\delta(p)
\mcW(k)
\widehat{J}(p,k)^*
=\frac{1}{2}\int_\T dk\; \mcW(k)\widehat{J}(0,k)^* ,
\end{split}\end{equation*}
where $\mcW(k)$  is  the Fourier transform of the correlation
function $\big<\psi(0)^*\psi(z)\big>$,
\begin{equation}\label{def:W}
\mcW(k)=\sum_{z\in\Z} e^{-2\pi i  k z}\big<\psi(0)^*\psi(z)\big> .
\end{equation}
By condition (d3), $\mcW(k)$ is well defined and in $L_1(\T)$.
Moreover, by translation invariance, $\mcW(k)$ is a real positive
function.

If we consider the deterministic dynamics only, then
$\mcW$ is preserved by the dynamics, i.e. $\partial_t \mcW(t)=0$.
This follows from eq. (\ref{evol}) for $\gamma=0$. Such property is
no longer true if the system evolves according to the full dynamics,
defined trough the generator $L=A+\ve\gamma S$. In order to observe
an effective change of the covariance, hence of the function $\mcW$, we
have to consider the time scale of order $\ve^{-1}$.
Denoting by $\mcW^\ve(t)=\mcW(t/\ve)$,
we obtain the following evolution equation,
\begin{equation}
\label{eq:evolW}
\partial_t\mcW^{\ve}(k,t)= \gamma (C\mcW^\ve)(k,t)
-\frac{\gamma}{2}[C(\mathcal Y^\ve+{\mathcal Y^\ve}^*)](k,t) ,
\end{equation}
where $C$ is the collision operator defined in $(\ref{C1})$, while
$\mathcal Y^\ve(k,t)$ is the Fourier transform of the correlation function
$\big<\psi(0)\psi(z)\big>$
at the rescaled time $t/\ve$,
$$
\mathcal Y^\ve(k,t)= \sum_{z\in\Z} e^{-2\pi i k
z}\big<\psi(t/\ve,0)\psi(t/\ve,z)\big> .
$$
As before, we prove that in the limit $\ve\to 0$ one obtains a closed
equation for $\mcW^\ve(t)$, as stated in the next theorem.

\begin{theorem}\label{theo:HBE}
Assume that the initial state satisfies the above conditions and
that $\mcW^\ve(k,0)= \mcW_0(k)$ is continuous on $\T$. Then,
$\forall k\in\T$, $t\in[0,T]$,
\begin{equation*}
\lim_{\ve\to 0}\mcW^\ve(k,t)= \mcW(k,t) ,
\end{equation*}
where $\mcW(k,t)$ satisfies the homogeneous Boltzmann equation
\begin{equation}\label{hbe}\begin{split}
\partial_t \mcW(k,t) & =\gamma (C\mcW)(k,t) , \\
\mcW(k,0) & = \mcW_0(k)
\end{split}\end{equation}
with $C$ defined in (\ref{C1}).
\end{theorem}
The proof of this theorem will be given in section \ref{sec:6.5}.


\subsection{Equilibrium time correlations}
\label{sec:equil-time-corr}

We consider the system in equilibrium and we denote by $\big< \cdot
\big>_T$ the average at respect to the equilibrium measure with
temperature $T$. This is a translation invariant Gaussian centered
measure with zero mean,  uniquely characterised through its
covariance
\begin{equation}
\label{meas_beta}
\begin{split}
\mcW(k) = T, \qquad \big<\psi(y)\psi(0)\big>_T
=0,\hspace{1cm}\forall y\in\Z .
\end{split}
\end{equation}
$\big< \cdot\big>_T$ is a stationary measure for the SDE (\ref{evol}).

Consider a function $g\in\ell_1(\Z)$ antisymmetric,  $g(z) = -
g(-z)$, and such that $\|\hat{g}/\omega\|_\infty<\infty$ and define
the function
\begin{equation*}
 \Phi =  \sum_{x\in\Z} g(x) p_x q_0.
\end{equation*}
The total time covariance is defined as
\begin{equation}
 \label{eq:stcf}
  \mathcal{F}^\ve(t)=\sum_{z\in\Z} \big<\Phi(t/\ve)\tau_z\Phi(0)\big>_T.
\end{equation}
We want to compute $\mathcal{F}^\ve(t)$ in the kinetic limit $\ve\to
0$.

Consider the centered translation invariant Gaussian measure defined
by the following covariance,
\begin{equation}\label{corr_beta-tau}
\mathcal{W}^{(T,\tau)}(k)=\frac{\omega(k)}{T^{-1}\omega(k)+i\tau
 \hat{g}(k)}\,,\hspace{1cm}
\mathcal{Y}^{(T,\tau)}(k)=0
\end{equation}
with $\tau >0$. Observe that $\mathcal{W}^{(\beta,\tau)}$ is a real,
continuous function which is positive for $\tau$ small enough.
Formally (\ref{corr_beta-tau}) corresponds to the perturbed measure
$$
Z^{-1}\exp\big[- T^{-1}H +\tau\sum_{z\in\Z}\tau_z\Phi\big].
$$
We denote by $\big<\cdot\big>_{(T,\tau)}$ its expectation.
\begin{lemma}For every $\ve>0$,
\begin{equation}\label{corr0}
\sum_{z\in\Z} \big<\Phi(t/\ve)\tau_z\Phi(0)\big>_T=
\lim_{\tau\to 0}\big<\Phi(t/\ve)\big>_{(T,\tau)}.
\end{equation}
\end{lemma}
The proof of this lemma is given in section \ref{proof-current} below.

By direct computation, at the rescaled time $t/\ve$,
\begin{equation*}\begin{split}
\big<\Phi(t/\ve)\big>_{(T,\tau)}=
-\frac{1}{2}\sum_{z\in\Z} g(z)\big<q_z p_0 - q_0 p_z\big>_{(T,\tau)}\\
=-i\int_\T dk\; \frac{\hat{g}(k)}{\omega(k)}\mathcal{W}^{(T,\tau)}(k,t/\ve).
\end{split}\end{equation*}
By Theorem \ref{theo:HBE}, $\mcW^{(T,\tau)}(k,t/\ve)\to \mcW(k,t)$
for $\ve\to 0$, where $\mcW(k,t)$ satisfies the homogeneous
Boltzmann equation (\ref{hbe}) with initial condition
$\mcW(k,0)=\omega(k)(T^{-1}\omega(k)+i\tau
 \hat{g}(k))^{-1}$.
It is easy to verify that for any bounded antisymmetric function $f$
on $\T$
\begin{equation*}
\int_\T dk\;(Cf)(k) \mcW(k,t)=-\int_\T dk\;\phi(k) \mcW(k,t),
\end{equation*}
with $\phi=-\hat\beta(k)=\frac{4}{3}\sin^2(\pi k)[1+2\cos^2(\pi
k)]$, see (\ref{hbeta}). Then
\begin{equation*}
-i\int_\T dk\;
\frac{\hat{g}(k)}{\omega(k)}\mathcal{W}^{(T,\tau)}(k,t)=
-i\int_\T dk\; \frac{\hat{g}(k)}{\omega(k)}
\mathcal{W}^{(T,\tau)}(k,0)e^{-\gamma\phi(k)t}
\end{equation*}
and finally, for $t\geq 0$,
\begin{equation}\begin{split}
\lim_{\ve\to 0}\mathcal{F}^\ve(t)&=\lim_{\tau\to 0}-\frac{i}{\tau}
\int_\T dk\;
\frac{\hat{g}(k)}{\omega(k)}\mathcal{W}^{(T,\tau)}(k,0)e^{-\gamma\phi(k)t}\\
&=T^2\int_\T dk\; \frac{|\hat{g}(k)|^2}{\omega(k)^2}e^{-\gamma\phi(k)t}.
\end{split}\end{equation}


\subsection{Energy current time correlation.}
The Hamiltonian energy current $\mathcal J$ is implicitly defined through the
conservation law
$$
A e_x=\tau_{x-1}\mathcal J - \tau_x \mathcal J.
$$
By direct computation $\mathcal J=\sum_{z>0}j_{0,z}$ with
$$j_{0,z}=-\frac{1}{2}\alpha(z)\sum_{y=0}^{z-1}(q_{z-y}p_{-y}-q_{-y}p_{z-y}).$$
Denoting by $\big<\cdot\big>$ the expectation value with respect to
some translation invariant centered Gaussian measure with covariance
$\mathcal W$, it is easy to see
that
\begin{equation*}\begin{split}
\big< \mathcal J\big>=-\frac{1}{2}\sum_{z>0}z \;\alpha(z)  \big<q_z p_0 - q_0 p_z
\big>\\
=\frac{1}{4\pi}\int_\T dk\;
\frac{\hat{\alpha}'(k)}{\omega(k)}\mathcal{W}(k)\\
=\frac{1}{2\pi}\int_\T dk\;\omega'(k)\mathcal{W}(k).
\end{split}\end{equation*}
Let us denote by $\mathcal C^\ve (t)$ the energy time correlation
function on the kinetic time scale $t/\ve$ at temperature $T$,
$$
\mathcal C^\ve (t)=\sum_{x\in\Z} \big<\mathcal J(t/\ve)\tau_x \mathcal
J(0) \big>_T. 
$$
Using the translation invariance of the Gaussian measure
$\big<\cdot\big>_T$ we have
\begin{equation*}\begin{split}
\mathcal C^\ve (t)&=\sum_{x\in\Z} \sum_{z\in\Z}
\frac{z}{4}\alpha(z)\sum_{z'\in\Z}
\frac{z'}{4}\alpha(z')\big<(q_zp_0-q_0p_z)(t/\ve)\\
&\hspace{20pt}(q_{x+z'}p_x-q_xp_{x+z'})(0)
\big>_T
=\sum_{x\in\Z}  \big<\tilde{\mathcal J}(t/\ve)\tau_x \tilde{\mathcal
  J}(0) \big>_T, 
\end{split}\end{equation*}
where
$$
\tilde{\mathcal J}=-\frac{1}{4}\sum_{z\in\Z} z\; \alpha(z)(q_zp_0-q_0p_z).
$$
Using the results of the previous subsection we arrive at
\begin{equation}\label{67}
\begin{split}
\lim_{\ve\to 0}\mathcal{C}^\ve(t)&=\frac{T^{2}}{(4\pi)^2}\int_\T
dk\;\frac{|\hat{\alpha}'(k)|^2}{\omega(k)^2}
e^{-\gamma\phi(k)|t|}\\
&=\frac{T^{2}}{4\pi^2}\int_\T dk\;|\omega'(k)|^2e^{-\gamma\phi(k)|t|}
\end{split}\end{equation}
for all $t\in\R$.

\bigskip

\subsection{Thermal Conductivity}
\label{sec:thermal-conductivity}

A standard definition of the thermal conductivity,
$\kappa^{(\varepsilon)}$, is by means of the Green-Kubo formula as
\begin{equation}\label{68}
\kappa^{(\varepsilon)}= \frac{1}{T^2}\int^\infty_0 dt\;
\mathcal{C}^\varepsilon(\varepsilon t) .
\end{equation}
We refer to \cite{bborev,bo} for details. In general, the kinetic
limit provides the lowest order approximation in $\varepsilon$ as
\begin{equation*}
\kappa^{(\varepsilon)}=\varepsilon^{-1}\kappa^{(0)}+\mathcal{O}(1).
\end{equation*}
Inserting in (\ref{68}) the limit (\ref{67}) thus yields
\begin{equation}\label{69}
\kappa^{(0)}=\frac{1}{4\pi^2}\int_\pi dk
\frac{|\omega'(k)|^2}{\gamma\phi(k)}.
\end{equation}
For the pinned case $\kappa^{(0)}<\infty$. 
$\kappa^{(\varepsilon)}$ has been computed in  \cite{bborev} with the result
\begin{equation*}
\kappa^{(\varepsilon)}=\varepsilon^{-1}\kappa^{(0)}+\varepsilon\gamma.
\end{equation*}
Thus, somewhat unexpectedly, the kinetic theory captures already the
main details of the conductivity.

For the unpinned case $\kappa^{(0)}=\infty$, hence
$\kappa^{(\varepsilon)}=\infty$, for $d=1,2$. The Boltzmann equation
(\ref{be}) provides a simple explanation. For small $k$,
$\omega(k)=2\pi c |k|$. Thus small $k$ phonons travel with speed
$c$. On the other hand, the collision rate vanishes as $k^2$ for
small $k$, see (\ref{hbeta}). Thus at small $k$ there are only very
few collisions which, together with $c>0$, is responsible for the
divergent conductivity. The positional part $X(t)$ of the process
consists mostly of very long stretches of uniform motion. In fact on
a large scale $X(t)$ is governed by a symmetric Levy process of
index $\alpha=3/2$, see \cite{JO} for details.


\subsection{Proof of theorem \ref{theo:HBE}}\label{sec:6.5}
The proof of the theorem \ref{theo:HBE} is analogous to the proof
above. We only have to control that $\mcW^\ve(t)$ and $\mcY^\ve(t)$
are well defined for every $t\in[0,T]$. This is stated in the next
lemma.
\begin{lemma}
Let the  conditions (d1-d3) hold.
Then $\mcW^\ve(t), \mcY^\ve(t)\in L_1(\T)$ for every $t\in[0,T]$.
\end{lemma}
\emph{Proof.} By similar computations as above we find that
$\mcW^\ve(t), \mcY^\ve(t)$ satisfy the following evolution equations,
\begin{equation*}\begin{split}
&\partial_t \mcW^\ve(k,t)=  \gamma
(C\mcW^\ve)(k,t)-\frac{\gamma}{2}(C(\mcY^\ve+{\mcY^\ve}^*))(k,t),\\
&\partial_t \mcY^\ve(k,t)=  -\frac{2i\omega(k)}{\ve}\mcY^\ve(k,t)
+\frac{\gamma}{2}\hat{\beta}(k)(\mcY^\ve+{\mcY^\ve}^*)(k,t)\\
&\hspace{60pt} +\frac{\gamma}{2}(C(\mcY^\ve+{\mcY^\ve}^*))(k,t)
-\frac{\gamma}{2}C(\mcW^\ve(k,t)+\mcW^\ve(-k,t)).
\end{split}\end{equation*}
In particular by Duhamel's formula we can rewrite the second equation
as
\begin{equation*}\begin{split}
\mcY^\ve(k,t)= & \gamma \int_0^t ds\; e^{
 -2i\omega(k)(t-s)/\ve}\big(
\frac{1}{2}\hat{\beta}(k)(\mcY^\ve+{\mcY^\ve}^*)(k,s)\\
& +\frac{1}{2}(C(\mcY^\ve+{\mcY^\ve}^*))(k,s)
-\frac{1}{2}C(\mcW^\ve(k,s)+\mcW^\ve(-k,s))
\big).
\end{split}\end{equation*}
Then we get the following bounds
\begin{equation*}\begin{split}
|\mcW^\ve(k,t)| & \leq \mcW(k) + \gamma c_1\int_0^t ds\;
\big[|\mcW^\ve(k,s)|+|\mcY^\ve(k,s)|\big]\\
& + \gamma c_2\int_0^t ds \int_\T dk\;
\big(|\mcW^\ve(k,s)|+|\mcY^\ve(k,s)|\big)\\
|\mcY^\ve(k,t)| & \leq + \gamma c_3\int_0^t ds\;
\big[|\mcW^\ve(k,s)|+|\mcW^\ve(-k,s)|+|\mcY^\ve(k,s)|\big]\\
& + \gamma c_4\int_0^t ds \int_\T dk\;
\big(|\mcW^\ve(k,s)|+|\mcY^\ve(k,s)|\big)\\
\end{split}\end{equation*}
and finally
\begin{equation*}\begin{split}
\int_\T dk\;
\big[|\mcW^\ve(k,t)|+|\mcY^\ve(k,t)|\big]\leq &\int_\T dk\;\mcW(k)+\\
& \gamma c_5\int_0^tds \int_\T dk\;
\big(|\mcW^\ve(k,s)|+|\mcY^\ve(k,s)|\big).
\end{split}\end{equation*}
By Gronwall's lemma
$$
\int_\T dk\; \big[|\mcW^\ve(k,t)|+|\mcY^\ve(k,t)|\big]\leq e^{\gamma
c_5 t} \int_\T dk\;\mcW(k) .
$$
\qed

\medskip

\subsection{Proof of (\ref{corr0})}\label{proof-current}
We define the generator of the speeded up process as
$$L_\ve=\ve^{-1}L=\ve^{-1}A +S.$$
Let us denote  by $\big< \cdot \big>_T$ the average with respect to
the Gaussian measure with zero mean and covariance
(\ref{meas_beta}), and by $\big< \cdot \big>_{(\beta,\tau)}$ the
Gaussian measure with zero mean and covariance
(\ref{corr_beta-tau}). We consider the Laplace transforms of
$\big<\Phi(t/\ve)\big>_{(T,\tau)}$ and
 $\sum_{x\in\mathbb{Z}}\big<\Phi(t/\ve)\tau_x\Phi(0)\big>_T$,
\begin{equation}\begin{split}
&\int_0^\infty dt\;e^{-\lambda t}\big<\Phi(t/\ve)\big>_{(T,\tau)}= 
\big<(\lambda-L_\ve)^{-1}\Phi\big>_{(T,\tau)}=\big<u_\lambda\big>_{(T,\tau)},\\
&\int_0^\infty dt\; e^{-\lambda t}\sum_x\big<\Phi(t/\ve)\tau_x \Phi(0)\big>_T= 
\sum_{x\in\Z} \big<[(\lambda-L_\ve)^{-1}\Phi]\tau_x\Phi\big>_T\\&\hspace{80pt} = 
\sum_{x\in\Z} \big<u_\lambda\tau_x\Phi\big>_T ,
\end{split}\end{equation}
where
$u_\lambda=\sum_{z\in\Z} f_\lambda(z)\; q_0 p_z$
with $f_\lambda$ the solution of the equation
$$
\lambda\; f_\lambda(z)- \frac{\gamma}{6}\Delta(4 f_\lambda(z) +f_\lambda(z+1)
+ f_\lambda(z-1)) = g(z).
$$
Observe that $u_\lambda$ does not depend on $\ve$, because $L_\ve
u_\lambda=\ve^{-1}\nabla F+\gamma S u_\lambda$ with  $F$ some
non-local function. Since $f_\lambda$ is antisymmetric, and by
translation invariance of the measure, the gradient term does not
contribute. We have
\begin{equation*}\begin{split}
\big<u_\lambda\big>_{T,\tau}&=\sum_{z\in\Z} f_\lambda(z)\big<q_0
p_z\big>_{T,\tau}=
-\frac{1}{2}\sum_{z\in\Z} f_\lambda(z)\big<q_z
p_0-q_0
p_z\big>_{T,\tau}\\
& =-i\int_\T dk \; \hat{f}_\lambda(k)\frac{1}{\omega(k)}
 \mcW ^{(T,\tau)}(k)\\ & =
-i \int_\T dk\; \frac{\hat{g}(k)}{\lambda+\gamma\phi(k)}
\frac{1}{T^{-1}\omega(k)+i\tau\hat{g}(k)}
\end{split}\end{equation*}
with $\phi(k)=-\hat{\beta}(k)$. For every positive $\lambda$, the
right hand side is finite for every $\tau$. In particular
$$
\lim_{\tau\to
 0}\frac{1}{\tau}\big<u_\lambda\big>_{T,\tau}= T^{2}\int_\T dk\;
\frac{|\hat{g}(k)|^2}{\omega(k)^2(\lambda+\gamma\phi(k))}\,.
$$
In the same way
\begin{equation*}\begin{split}
&\sum_{x\in\Z}
\big<u_\lambda\;\tau_x\Phi\big>_T= T\sum_{y\in\Z}\sum_{z\in\Z}
f_\lambda(y) g(z)\big< q_0 q_{y-z})\big>_T\\
&\hspace{10pt}= T\int_\T dk\;
\hat{f}_\lambda(k)^*\frac{\hat{g}(k)}{\omega(k)^2}\mcW
^{(\beta)}(k) = T^{2} \int_\T dk\;
\frac{|\hat{g}(k)|^2}{\omega(k)^2(\lambda+\phi(k))}\,,
\end{split}\end{equation*}
where Parseval's  identity is used.
\qed


\section{Appendix}
\label{app}


\begin{proposition}\label{prop:bound}
Under the assumption
(b3),
\begin{equation}
\big|\big<J,W^\ve \big>\big|\leq C_J 
\end{equation}
for every test function
$J\in \mathcal S (\R^d\times\T^d)$ and for every $\ve>0$.
\end{proposition}
\emph{Proof}. For any test function
$J\in\mathcal{S}(\R^d\times\T^d)$ we denote by $\|J\|_{n,\infty}$
the following norm,
$$
\|J\|_{n,\infty}=\sup_{i,j\in\{1,\ldots,d\}}\sup_{\br\in\R^d,
\bk\in\T^d} \big|\sum_{i_1=1}^d\ldots\sum_{i_n=1}^d 
\partial_{k_{i_1}}\ldots\partial_{k_{i_n}}J_{i,j}(\br,\bk)
\big|.
$$
Let us  define
$$
\mathcal{W}^\ve[J](\by,\by',i,j)=\int_{\T^d} d\bk \;e^{2\pi i
\bk\cdot(\by'-\by)} J_{j,i}\ve(\by'+\by)/2,\bk)^*,
$$
for every $\by,\by'\in\Z^d$, $i,j=1,\ldots,d$. Integrating by part
in $\bk$ for $(d+1)$ times, we get
\begin{equation*}
|\mathcal{W}^\ve[J](\by,\by',i,j)|\leq \frac{1}{(2\pi)^{d+1}}
\frac{1}{|(y'_1-y_1)+\ldots+(y'_d-y_d)|^{d+1}}\|J\|_{d+1,\infty},
\end{equation*}
where $y_i$ denotes the $i$-th component of the vector $\by$. We denote
by
$$\big<\by'-\by\big>^{d+1}=|(y'_1-y_1)+\ldots+(y'_d-y_d)|^{d+1}.$$
By Schwarz inequality
\begin{equation*}\begin{split}
&\big|\big<J,W^\ve \big>\big|\\&\hspace{10pt}\leq
(\ve/2)^d\sum_{\by,\by'\in\Z^d}\sum_{i,j=1}^{d} 
[\big<|\psi_i(\by)|^2\big>_\ve]^{1/2}[\big<|\psi_j(\by')|^2\big>_\ve]^{1/2}
|\mathcal{W}^\ve[J](\by,\by',i,j)|\\
&\hspace{10pt}\leq (\ve/2)^d\big<\|\bpsi\|^2\big>_\ve
\sum_{\bz\in\Z^d}\frac{1}{\big<\bz\big>^{d+1}}
c_0 \|J\|_{d+1,\infty} \leq c \|J\|_{d+1,\infty},
\end{split}\end{equation*}
where in the last inequality we used $\ve^d\big<\|\bpsi\|^2\big>\leq
K$ with $K$ positive.
\qed

\bigskip

\begin{lemma}\label{lemma_omega}
Let Assumption \ref{alpha_assump} hold with $\hat{\alpha}(0)=0$
$\mathrm{(unpinned}$ $\mathrm{case)}$. The following assertions hold.
\begin{itemize}
\item[(i)] There are constants $C_1$, $C_2$, $C_3$ such that $\forall \bk\in\T^d$
\begin{equation}
|\nabla\hat{\alpha}(\bk)|\leq C_1|\bk|,
\hspace{1cm}C_2|\bk|\leq\omega(\bk)\leq C_3|\bk|.
\end{equation}
In addition,
$\|\nabla\omega\|_\infty<\infty$.
\item[(ii)] For all $\bk\in\T^d$  and
$\bp\in\R^d$, there is a positive constant $C$ such that
\begin{equation}
\ve^{-1}\big|\omega(\bk+\ve\bp/2)-\omega(\bk-\ve\bp/2)\big|\leq C|\bp|
\end{equation}
for every $\ve>0$.
\item[(iii)] For $\bp\in\R^d$, there is a positive constant $C_4$ such
that for all $\bk\in\T^d$, $|\bk|>\ve|\bp|$, with $\ve>0$ one has
\begin{equation}\label{iii}
\big|\ve^{-1}[\omega(\bk+\ve\bp/2)-\omega(\bk-\ve\bp/2)]-\bp\cdot\nabla\omega(\bk)\big|\leq
\ve C_4\frac{|\bp|^2}{|\bk|}\,.
\end{equation}
\end{itemize}
\end{lemma}

\emph{Proof.} The first inequality of item (i) follows from a Taylor
expansion of $\hat{\alpha}$ around zero, using the fact that
$\nabla\hat{\alpha}(0)=0$ and $\|D^2\hat{\alpha}\|_\infty<\infty$.
Using the same argument we have $\omega(\bk)\leq C_3|\bk|$, since
$\omega(\bk)=\hat{\alpha}(\bk)^{1/2}$ and $\hat{\alpha}(\bk)\leq
C|\bk|^2$, with $C>0$. Let us denote with $A_0$ the Hessian of
$\hat{\alpha}$ at $\bk =0$. By assumption (a4), $\bk\cdot
A_0\bk>c|\bk|^2$. Moreover, there is a $\delta>0$ such that for
every $|\bk|<\delta$, $(|\hat{\alpha}(\bk)-\frac{1}{2}\bk\cdot
A_0\bk|)^2<\frac{1}{4}\bk\cdot A_0\bk$. Then
$\omega(\bk)=((\hat{\alpha}(\bk)-\frac{1}{2}\bk\cdot A_0\bk)
+\frac{1}{2} \bk\cdot A_0\bk)^{1/2}\geq \frac{c}{2}|\bk|^2$ if $
|\bk|<\delta$. For $|\bk|\geq \delta$, $\omega$ is strictly positive
and there is a constant $c'$ such that $\omega(\bk)\geq c'|\bk|$ if
$|\bk|\geq \delta$, which proves the second inequality. To prove the
last one it is enough to observe that for all $\bk\neq 0$,
$$|\nabla\omega(\bk)|=\frac{1}{2}\frac{|\nabla\hat{\alpha}(\bk)|}{\omega(\bk)}\leq
\frac{1}{2}\frac{C_1}{C_2}<\infty.$$
Hence $\|\omega\|_\infty<\infty$.

Let us prove item (iii). We observe that the function on the left hand side of
(\ref{iii}) is zero if $|\bp|=0$ and  we have to consider only the
case $|\bp|>0$. Since we are assuming $|\bk|>\ve|\bp|$, it follows
that we need to discuss the case $|\bk|>0$. Observe that for any
$s\in(0,\ve]$, $|\bk\pm \frac{1}{2}s\bp|\geq \frac{1}{2}|\bk|>0$ if
$|\bk|>\ve|\bp|$ and the function
$\omega(\bk+\frac{1}{2}\ve\bp)-\omega(\bk-\frac{1}{2}\ve\bp)$ is
$C^\infty$ in this range. In particular
\begin{equation*}\begin{split}
&\ve^{-1}\omega(\bk+\ve\bp/2)-\omega(\bk-
\ve\bp/2)-\bp\cdot\nabla\omega(\bk)\\
&\hspace{10pt}=\ve \left(\frac{\bp}{2}\cdot\nabla\right)^2\omega(\bk+s\bp/2)-
\ve \left(\frac{\bp}{2}\cdot\nabla\right)^2\omega(\bk-\tilde{s}\bp/2)
\end{split}\end{equation*}
with $s,\tilde{s}\in(0,\ve)$,
where, denoting  $\bk_\pm=\bk\pm \frac{1}{2}s\bp$,
\begin{equation*}
\left(\bp\cdot\nabla\right)^2\omega(\bk_\pm)=\frac 12
\frac{1}{\omega(\bk_\pm)}
\left(\bp\cdot\nabla\right)^2\hat{\alpha}(\bk_\pm)-\frac{1}{4}
\frac{1}{\omega(\bk_\pm)^3}
\left(\bp\cdot\nabla\hat{\alpha}(\bk_\pm)\right)^2
\end{equation*}
and, using item (i),
\begin{equation*}
\left|\left(\bp\cdot\nabla\right)^2\omega(\bk_\pm) \right|\leq 4
C_4\frac{|\bp|^2}{|\bk|}\,.
\end{equation*}
This prove item (iii).

Item (ii) follows from (iii) if $|\bk|>\ve|\bp|$. If $|\bk|\leq
\ve|\bp|$ we use the bound
\begin{equation*}
\ve^{-1}|\omega(\bk+\ve\bp/2)-\omega(\bk-\ve\bp/2)|\leq
C_3\ve^{-1}(|\bk+\ve\bp/2|+|\bk-\ve\bp/2|)\leq C|\bp|.
\end{equation*}


\subsection{Proof of (\ref{eq:Ypsi}).}\label{proof:SW}

First of all we observe that
\begin{equation}\label{syy}
\begin{split}
&\sum_{z\in\Z} [Y_z\psi(y')^*][Y_z\psi(y)]
=  [Y_{y+1}\psi(y')^*][Y_{y+1}\psi(y)] \\
&\hspace{60pt}+
[Y_y\psi(y')^*][Y_y\psi(y)]
+[Y_{y-1}\psi(y')^*][Y_{y-1}\psi(y)],
\end{split}\end{equation}
where
\begin{equation*}\begin{split}
&Y_{y+1}\psi(y)=\frac{i}{\sqrt{2}}(p_{y+1}-p_{y+2}),\quad
Y_y\psi(y)=\frac{i}{\sqrt{2}}(p_{y+1}-p_{y-1}),\\
&\hspace{0pt}Y_{y-1}\psi(y)=\frac{i}{\sqrt{2}}(p_{y-2}-p_{y-1}).
\hspace{1cm}
\end{split}\end{equation*}
Then (\ref{syy}) is equal to
\begin{equation}\label{syy2}
\begin{split}
&(p_{y+1} - p_{y+2})\big[ (p_{y+1} - p_{y+2}) \delta_{y',y}
+ (p_{y+2} - p_{y}) \delta_{y',y+1}\\ 
&+ (p_{y} - p_{y+1}) \delta_{y',y+2}
\big] 
+   (p_{y+1} - p_{y-1})\big[ (p_{y} - p_{y+1}) \delta_{y',y-1}
\\&+ (p_{y+1} - p_{y-1}) \delta_{y',y} + 
(p_{y-1} - p_{y}) \delta_{y',y+1}
\big] 
+ (p_{y-2} - p_{y-1})\\
&\times\big[ (p_{y-1} - p_{y}) \delta_{y',y-2}
+ (p_{y} - p_{y-2}) \delta_{y',y-1} + (p_{y-2} - p_{y-1}) \delta_{y',y}
\big] \\
&=   (p_{y+1} - p_{y+2}) (p_{y} - p_{y+1}) \delta_{y',y+2}
+  (p_{y-2} - p_{y-1}) (p_{y-1} - p_{y}) \delta_{y',y-2}\\
&+   \left[ (p_{y+1} - p_{y+2}) (p_{y+2} - p_{y})
+(p_{y+1} - p_{y-1}) (p_{y-1} - p_{y}) \right] \delta_{y',y+1} \\
&+   \left[ (p_{y+1} - p_{y-1}) (p_{y} - p_{y+1})
+(p_{y-2} - p_{y-1}) (p_{y} - p_{y-2}) \right] \delta_{y',y-1} \\
&+    \left[ (p_{y+1} - p_{y+2})^2
+(p_{y+1} - p_{y-1})^2 + (p_{y-2} - p_{y-1})^2 \right] \delta_{y',y} \\
&=  2\sum_{r=-2}^2 A_{y,r} \delta_{y',y+r}\\
\end{split}\end{equation}
and we can write
\begin{equation*}
\begin{split}
&\frac {1}{3}(\ve/2) \sum_{y,y',z\in\Z} \big<
[Y_z\psi(y')^*][Y_z\psi(y)]\big>_\ve \tilde
J(\ve (y'+y)/2, y - y')^* \\
&\hspace{10pt}= \frac {1}{3}(\ve/2) \sum_{y\in\Z}\sum_{r=-2}^2 A_{y,r} \tilde
J(\ve (y+ r/2) , -r)^*.
\end{split}
\end{equation*}
Expanding the $A_{y,r}$ and identifying terms that differ by
translations, we arrive at
\begin{equation*}
\begin{split}
&\frac {1}{3}(\ve/2) \sum_{y,y',z\in\Z} \big<
[Y_z\psi(y')^*][Y_z\psi(y)]\big>_\ve \tilde J(\ve (y'+y)/2, y -
y')^* \\ &\hspace{10pt} = \frac {\ve}{6} \sum_{y\in\Z} \Big[ \big<2p_y p_{y+1} -
p_{y}p_{y+2} - p_y^2\big>_\ve (\tilde J(\ve y , 2)^* +\tilde J(\ve
y, -2)^*)\\ 
&\hspace{20pt}  + \big<2p_{y}p_{y+2} - 2p_y^2\big>_\ve (\tilde
J(\ve y , 1)^* + \tilde J(\ve y , -1)^*)\\ 
&\hspace{20pt}  + \big<-4p_y p_{y+1}
- 2p_{y}p_{y+2} + 6p_y^2\big>_\ve \tilde J(\ve y , 0)^*\Big] +
\mathcal{O}(\ve),
\end{split}
\end{equation*}
where we have used the smoothness of $\tilde J$ in $x\in\R$.
We can rewrite the last expression as
\begin{equation*}
(\ve/2) \sum_{y\in\Z} \sum_{z,u=-2}^2 \alpha(z,u) \big<p_y p_{y+z}\big>_\ve
\tilde J( \ve y , u)^* + \mathcal{O}(\ve),
\end{equation*}
where $ \alpha(z,u) = \alpha(z,-u) = \alpha(u,z)=\alpha(-u,z)$ and is given by
\begin{equation}\label{eq:alpha}
\begin{array}{lll}
\alpha (0,0) = 1, &\alpha (0,1)=-1/3, &\alpha(0,2)=-1/6\\
&\alpha(1,1)=0, &\alpha(1,2)=1/6\\
& &\alpha (2,2) = -1/12.
\end{array}
\end{equation}


\bigskip\bigskip
\end{document}